\documentclass{svmult}
\usepackage{amssymb,amsmath,amsfonts,latexsym}
\usepackage{mathrsfs}
\usepackage{cite, hyperref}
\newtheorem{Question}{Question}
\newcommand{\B}{\mathcal{B}}
\newcommand{\C}{ \mathbb{C}}
\renewcommand{\D}{ \mathbb{D}}
\newcommand{\dD}{ \partial\mathbb{D}}
\newcommand{\Z}{ \mathbb{Z}}

\newcommand{\ran}{\operatorname{ran}}
\renewcommand{\Re}{\operatorname{Re}}

\newcommand{\norm}[1]{\| #1 \|}
\newcommand{\inner}[1]{\langle #1 \rangle}
\newcommand{\A}{\mathcal{A}}
\renewcommand{\E}{\mathcal{E}}

\newcommand{\M}{\mathcal{M}}
\newcommand{\N}{\mathbb{N}}

\newcommand{\h}{\mathcal{H}}
\newcommand{\K}{\mathcal{K}}

\renewcommand{\vec}[1]{{\bf #1}}

\renewcommand{\phi}{\varphi}

\newcommand{\TTO}{\mathscr{T}}

\allowdisplaybreaks

\begin{document}

    \title*{Recent Progress on Truncated Toeplitz Operators}

    \author{Stephan Ramon Garcia and William T. Ross}
    \institute{Stephan Ramon Garcia \at   Department of Mathematics, Pomona College, Claremont, CA 91711, USA \email{Stephan.Garcia@pomona.edu} \and
   William T. Ross \at   Department of Mathematics and Computer Science, University of Richmond,  Richmond, VA 23173, USA \email{wross@richmond.edu}}
   
    \maketitle

\abstract{
	This paper is a survey on the emerging theory of truncated Toeplitz operators.
	We begin with a brief introduction to the subject and then highlight the many recent
	developments in the field since Sarason's seminal paper \cite{Sarason} from 2007.}

%\thanks{First author partially supported by National Science Foundation Grant DMS-1001614.}

\section{Introduction}

	Although the subject has deep classical roots, the systematic study of \emph{truncated Toeplitz operators}
	for their own sake was only recently spurred by the seminal 2007 paper of Sarason \cite{Sarason}.
	In part motivated by several of the problems posed in the aforementioned article,
	the area has undergone vigorous development during the past several years
	\cite{Sarason, BCFMT,BBK,MR2597679,TTOSIUES,MR2440673,VPSBF,CCO,ATTO,GRW,G-P,G-P-II,NLEPHS,GRNorm,
	RH,MR2418122,MR2468883,Sed,SedlockThesis,STZ}.  While several of the initial
	questions raised by Sarason have now been resolved, the study of truncated Toeplitz operators
	has nevertheless proven to be fertile ground, spawning both new questions and unexpected results.
	In this survey, we aim to keep the interested reader up to date and to give the uninitiated reader a historical
	overview and a summary of the important results and major developments in this area.
	
	Our survey commences in Section \ref{SectionPreliminaries} with an extensive treatment of the basic definitions,
	theorems, and techniques of the area.  Consequently, we shall be brief in this introduction and simply declare that
	a \emph{truncated Toeplitz operator} is the compression $A_{\phi}^u:\K_u\to\K_u$ of a classical
	Toeplitz operator $T_{\phi}$ to a shift coinvariant subspace $\K_u:=H^2\ominus uH^2$
	of the classical Hardy space $H^2$.  Here $u$ denotes a nonconstant inner function
	and we write $A^u_{\phi}f = P_u(\phi f)$ where $P_u$ denotes the orthogonal projection from $L^2$ onto $\K_u$.
	Interestingly, the study of potentially unbounded truncated Toeplitz operators, having symbols $\phi$
	in $L^2$ as opposed to $L^{\infty}$, has proven to be spectacularly fruitful.
	Indeed, a number of important questions in the area revolve around this theme
	(e.g., the results of Section \ref{SectionBoundedSymbol}).

	Before proceeding, let us first recall several instances
	where truncated Toeplitz operators have appeared in the literature.  This will not only
	provide a historical perspective on the subject, but it will also illustrate the fact that
	truncated Toeplitz operators, in various guises, form the foundations of much of modern
	function-related operator theory.
	
	Let us begin with the poweful Sz.-Nagy-Foia\c{s} model theory for Hilbert space contractions,
	the simplest manifestation of which is the \emph{compressed shift} $A_z^u$ \cite{Bercovici, N1,N2,N3,RR}.
	To be more specific, every Hilbert space contraction $T$ having defect indices $(1,1)$
	and such that
	$\lim_{n \to \infty} T{^*}^n = 0$ (SOT) is unitarily equivalent to $A_z^u$
	for some inner function $u$.  Natural generalizations of this result are available to treat
	contractions with arbitrary defect indices by employing the machinery of vector-valued model spaces
	and operator-valued inner functions.

	In his approach to the \emph{Gelfand problem}, that is the characterization of the invariant
	subspace lattice $\mbox{Lat} V$ of the Volterra integration operator
	\begin{equation}\label{eq-Volterra}
		[Vf](x) = \int_0^xf(y)\,dy
	\end{equation}
	on $L^2[0,1]$, Sarason noted that the Volterra operator
	is unitarily equivalent to the Cayley transform of the compressed shift  $A_z^u$
	corresponding to the atomic inner function $u(z) = \exp( \frac{z+1}{z-1})$ \cite{MR0192355}.
	This equivalence was then used in conjunction with Beurling's Theorem to demonstrate
	the unicellularity of $\mbox{Lat} V$ \cite{MR0358396, MR0192355,N1,N3}.
	Interestingly, it turns out that the Volterra operator, and truncated Toeplitz operators in general,
	are natural examples of \emph{complex symmetric operators},
	a large class of Hilbert space operators which has also undergone significant development in recent
	years \cite{Chevrot, CCO, G-P, G-P-II, ESCSO,VPSBF,ATCSO,ONCPCSO,OCCSO,
	SNCSO, CSPI, Gilbreath, JKLL,JKLL2, ZL, ZLZ, Zag,WXH,Wang, ZLJ}.  This link between
	truncated Toeplitz operators and complex symmetric operators is explored in
	Section \ref{SectionCSO}.
	
	Sarason himself identified the commutant of $A^{u}_{z}$ as the
	set $\{A^{u}_{\phi}: \phi \in H^{\infty}\}$ of all \emph{analytic} truncated Toeplitz operators.
	He also obtained an $H^{\infty}$ functional calculus for the compressed shift, establishing that
	$\phi(A_z^u) = A_{\phi}^u$ holds for all $\phi$ in $H^{\infty}$ \cite{Sarason-NF}.
	These seminal observations mark the beginning of the so-called \emph{commutant lifting} theory, which
	has been developed to great effect over the ensuing decades \cite{MR1120546,MR2158502, SarasonCL}.
	Moreover, these techniques have given new perspectives on several classical problems
	in complex function theory.  For instance, the Carath\'eodory and Pick problems lead one naturally
	to consider lower triangular Toeplitz matrices (i.e., analytic truncated Toeplitz operators on $\K_{z^n}$)
	and the backward shift on the span of a finite collection of Cauchy kernels (i.e., $A_{\overline{z}}^u$ on
	a finite dimensional model space $\K_u$).  We refer the reader to the text \cite{AMC} which treats these
	problems in greater detail and generality.

	Toeplitz matrices, which can be viewed as truncated Toeplitz operators on $\K_{z^n}$,
	have long been the subject of intense research.
	We make no attempt to give even a superficial discussion of this immense topic.
	Instead, we merely refer the reader to several recent texts which analyze various
	aspects of this fascinating subject.
	The pseudospectral properties of Toeplitz matrices are explored in \cite{Trefethen}.	
	The asymptotic analysis of Toeplitz operators on $H^2$ via large truncated Toeplitz matrices
	is the focus of \cite{BS}.
	The role played by Toeplitz determinants in the study of orthogonal polynomials
	is discussed in \cite{OPUC1, OPUC2} and its relationship to random matrix theory
	is examined in \cite{MR1957518,MR2166467}.  Finally, we should also remark
	that a special class of Toeplitz matrices, namely circulant matrices,
	are a crucial ingredient in many aspects of numerical computing \cite{Davis}.

	We must also say a few words about the appearance of truncated Toeplitz operators
	in applications to control theory and electrical engineering.  In such contexts, extremal problems posed over $H^{\infty}$
	often appear.  It is well-known that the solution to many such problems can be obtained by computing
	the norm of an associated Hankel operator \cite{MR901235, MR971878}.  However,
	it turns out that many questions about Hankel
	operators can be phrased in terms of analytic truncated Toeplitz operators and, moreover,  this link has long
	been exploited \cite[eq.~2.9]{Peller}.  Changing directions somewhat, we remark that the
	\emph{skew Toeplitz operators} arising in $H^{\infty}$ control theory are closely related to
	selfadjoint truncated Toeplitz operators \cite{MR945002, MR1639646}.

	Among other things, Sarason's recent article \cite{Sarason} is notable for opening the general study of
	truncated Toeplitz operators, beyond the traditional confines of the analytic ($\phi \in H^{\infty}$)
	and co-analytic $(\overline{\phi} \in H^{\infty}$) cases and the limitations of the case $u = z^N$ (i.e., Toeplitz matrices),
	all of which are evidently well-studied in the literature.  By permitting arbitrary symbols in $L^{\infty}$,
	and indeed under some circumstances in $L^2$,
	an immense array of new theorems, novel observations, and difficult problems emerges.
	It is our aim in this article to provide an overview of the ensuing developments, with an eye toward
	promoting further research.  In particular, we make an effort to highlight open problems and unresolved issues
	which we hope will spur further advances.

\section{Preliminaries}\label{SectionPreliminaries}
	In this section we gather together some of the standard results on model spaces and Aleksandrov-Clark
	measures which will be necessary for what follows.  Since most of this material is familiar to those who have studied
	Sarason's article \cite{Sarason}, the following presentation is somewhat terse.  Indeed, it serves primarily as a
	review of the standard notations and conventions of the field.

\subsection{Basic Notation} Let $\D$ be the open unit disk, $\partial \D$ the unit circle, $m = d \theta/2\pi$ \emph{normalized} Lebesgue measure on $\dD$, and $L^p := L^p(\dD, m)$ be the standard Lebesgue spaces on $\dD$. For $0 < p < \infty$ we  use $H^p$ to denote the classical Hardy spaces of $\D$ and $H^{\infty}$ to denote the bounded analytic functions on $\D$. As is standard, we regard $H^p$ as a closed subspace of $L^p$ by identifying each $f \in H^p$ with its $m$-almost everywhere $L^p$  boundary function
$$f(\zeta) := \lim_{r \to 1^{-}} f(r \zeta), \quad \mbox{$m$-a.e. $\zeta \in \dD$}.$$ In the Hilbert space setting $H^2$ (or $L^2$) we denote the norm as $\|\cdot\|$ and the usual integral inner product by $\langle \cdot, \cdot \rangle$. On the rare occasions when we need to discuss $L^p$ norms we  will use $\|\cdot\|_{p}$.
We let $\widehat{\C}$ denote the Riemann sphere $\C \cup \{\infty\}$ and, for a set $A \subseteq \C$, we let $A^{-}$ denote the closure of $A$. For a subset $V \subset L^p$, we let $\overline{V} :=\{\overline{f}: f \in V\}$.
We interpret the Cauchy integral formula
	$$
		f(\lambda) = \int_{\dD} \frac{f(\zeta)}{1 - \overline{\zeta} \lambda} dm(\zeta),
	$$
	valid for all $f$ in $H^2$, in the context of reproducing kernel Hilbert spaces by writing
	$f(\lambda) = \inner{ f, c_{\lambda}}$ where
	\begin{equation}\label{eq-CauchyKernel}
		c_{\lambda}(z) := \frac{1}{1 - \overline{\lambda} z}
	\end{equation}
	denotes the \emph{Cauchy kernel} (also called the \emph{Szeg\H{o} kernel}).
	A short computation now reveals that the orthogonal projection $P_+$ from $L^2$ onto $H^2$
	(i.e., the \emph{Riesz projection}) satisfies
	$$
		[P_{+} f](\lambda) = \inner{ f, c_{\lambda} }
	$$
	for all $f \in L^2$ and $\lambda \in \D$.
	
\subsection{Model Spaces}
	%In the following, we let $H^2$ denote the classical Hardy space on the open unit disk
	%$\D$, standard references for which are \cite{Duren, Garnett, Hoffman}.  Letting
	%$d m = \frac{d \theta}{2 \pi}$ denote normalized Lebesgue measure on the unit circle $\dD$,
	%we may regard $H^2$ as a closed subspace of $L^2 := L^2(\dD,m)$
	%by identifying each $f$ in $H^2$ with its a.e.~defined nontangential boundary function on $\dD$.	

	Let $S: H^2 \to H^2$ denote the unilateral shift
	\begin{equation}\label{eq-ShiftOperator}
		[S f](z) = z f(z),
	\end{equation}
	and recall that Beurling's Theorem asserts that the nonzero $S$-invariant subspaces
	of $H^2$ are those of the form $uH^2$ for some inner function $u$.  Letting
	\begin{equation}\label{eq-BackwardShift}
		[S^{*} f](z) = \frac{f(z) - f(0)}{z}
	\end{equation}
	denote the backward shift operator, it follows that
	the proper $S^*$-invariant subspaces of $H^2$ are precisely those of the form
	\begin{equation}\label{eq-ModelSpace}
		\K_{u} := H^2 \ominus u H^2.
	\end{equation}
	The subspace \eqref{eq-ModelSpace} is called the \emph{model space}
	corresponding to the inner function $u$, the terminology stemming from the important role that $\K_u$ plays in the
	model theory for Hilbert space contractions \cite[Part C]{N3}.
	
	Although they will play only a small role in what follows, we should also mention that
	the backward shift invariant subspaces of the Hardy spaces $H^p$ for $0 < p < \infty$ are also known.
	In particular, for $1 \leqslant p < \infty$ the proper backward shift invariant subspaces of $H^p$ are all of the form
	\begin{equation}\label{eq-pModel}
		\mathcal{K}_{u}^{p} := H^{p} \cap u \overline{H^{p}_{0}},
	\end{equation}
	where $H^p_0$ denotes the subspace of $H^p$ consisting of those $H^p$ functions which vanish at the origin and where the right-hand side of \eqref{eq-pModel} is to be understood in terms of boundary functions on $\dD$.
	For further details and information on the more difficult case $0 < p < 1$, we refer the reader to the text \cite{CR} and
	the original article \cite{MR566739} of Aleksandrov.  For $p =2$, we often suppress
	the exponent and simply write $\K_u$ in place of $\K_u^2$.
	
\subsection{Pseudocontinuations}
	Since the initial definition \eqref{eq-pModel} of $\K_u^p$ is somewhat indirect, one might
	hope for a more concrete description of the functions belonging to $\K_u^p$.
	A convenient function-theoretic characterization of $\K_u^p$ is provided by
	the following important result.

	\begin{theorem}[Douglas-Shapiro-Shields \cite{DSS}] \label{Thm-DSS}
		If $1 \leqslant p < \infty$, then $f$ belongs to $\K_u^p$ if and only if
		there exists a $G \in H^p(\widehat{\C} \backslash \D^{-})$ which vanishes at infinity\footnote{Equivalently, $G(1/z) \in H^{p}_{0}$.}  such that
		\begin{equation*}
			\lim_{r \to 1^{+}} G(r \zeta) = \lim_{r \to 1^{-}} \frac{f}{u}(r \zeta)
		\end{equation*}
		for almost every $\zeta$ on $\dD$.
	\end{theorem}

	The function $G$ in the above theorem is called a \emph{pseudocontinuation} of $f/u$
	to $\widehat{\C} \backslash \D^{-}$.  We refer the reader
	to the references \cite{MR629832, CR, DSS, RS} for further information about pseudocontinations
	and their properties.  An explicit function theoretic parametrization of the spaces
	$\K_u^p$ is discussed in detail in \cite{CCO}.
	
\subsection{Kernel Functions, Conjugation, and Angular Derivatives}

	Letting $P_u$ denote the orthogonal projection from $L^2$ onto $\K_u$, we see that
	\begin{equation} \label{eq-Pu-formula}
		[P_{u} f](\lambda) = \inner{ f, k_{\lambda}}
	\end{equation}
	for each $\lambda$ in $\D$.  Here
	\begin{equation} \label{eq-rk-defn}
		k_{\lambda}(z) := \frac{1 - \overline{u(\lambda)} u(z)}{1 - \overline{\lambda} z}
	\end{equation}
	denotes the \emph{reproducing kernel} for $\K_u$.  In particular, this family of functions
	has the property that $f(\lambda)= \inner{f, k_{\lambda}}$ for every $f \in \K_{u}$ and $\lambda \in \D$.
	
	Each model space $\K_u$ carries a natural \emph{conjugation} (an isometric,
	conjugate-linear involution) $C: \K_u \to \K_u$, defined in terms of boundary functions by
	\begin{equation}\label{eq-ModelConjugation}
		[C f](\zeta) := \overline{ f(\zeta)\zeta} u(\zeta).
	\end{equation}
	For notational convenience, we sometimes denote the conjugate $Cf$ of $f$ by $\tilde{f}$.
	More information about conjugations in general, along with specific properties of the map
	\eqref{eq-ModelConjugation} can be found in \cite{CCO}.  	
	For the moment, we simply mention that the so-called \emph{conjugate kernels}
	\begin{equation}\label{eq-ConjugateKernel}
		[Ck_{\lambda}](z) = \frac{ u(z) - u(\lambda)}{z - \lambda}
	\end{equation}
	will be important in what follows.  In particular, observe that each conjugate kernel
	is a difference quotient for the inner function $u$.  We therefore expect that derivatives
	will soon enter the picture.
	
	\begin{definition}
		For an inner function $u$ and a point $\zeta$ on $\dD$ we say that $u$ has an
		\emph{angular derivative in the sense of Carath\'eodory} (ADC)
		\emph{at $\zeta$} if the nontangential limits of $u$ and $u'$ exist at $\zeta$ and  $|u(\zeta)| = 1$.
	\end{definition}

	The following theorem provides several useful characterizations of ADCs.
	
	\begin{theorem}[Ahern-Clark \cite{AC70}]  \label{AC-paper}
		For an inner function $u = b_{\Lambda} s_{\mu}$, where $b_{\Lambda}$
		is a Blaschke product with zeros $\Lambda = \{\lambda_n\}_{n=1}^{\infty}$, repeated
		according to multiplicity,  $s_{\mu}$ is a singular inner function
		with corresponding singular measure $\mu$,  and $\zeta \in \dD$, the following are equivalent:
		\begin{itemize}\addtolength{\itemsep}{0.5\baselineskip}
			\item[(i)]\quad Every $f \in \mathcal{K}_{u}$ has a nontangential limit at $\zeta$.
			\item[(ii)]\quad For every $f \in \mathcal{K}_u$, $f(\lambda)$ is bounded as $\lambda \to \zeta$ nontangentially.
			\item[(iii)]\quad $u$ has an ADC at $\zeta$.
			\item[(iv)]\quad The function
				\begin{equation}\label{eq-BoundaryKernel}
					k_{\zeta}(z) = \frac{1 - \overline{u(\zeta)} u(z)}{1 - \overline{\zeta} z},
				\end{equation}
				belongs to $H^2$.
			\item[(v)]\quad The following condition holds:
				\begin{equation} \label{AC1-2}
					\sum_{n \geq 1} \frac{1 - |\lambda_{n}|^2}{|\zeta - \lambda_{n}|^2} +
					\int_{\dD} \frac{d \mu(\xi)}{|\xi - \zeta|^2} < \infty.
				\end{equation}
		\end{itemize}
	\end{theorem}

	In fact, the preceding is only a partial statement of the Ahern-Clark result, for there are several other
	additional conditions which can be appended to the preceding list.   Moreover, they also characterized the
	existence of nontangential boundary limits of the derivatives (up to a given order) of functions in $\mathcal{K}_u$.
	An extension of Theorem \ref{AC-paper} to the spaces $\mathcal{K}_u^{p}$ is due to Cohn \cite{Cohn}.
	
	Among other things, Theorem \ref{AC-paper} tells us that whenever $u$ has an ADC at a point
	$\zeta$ on $\dD$, then the functions \eqref{eq-BoundaryKernel},
	are reproducing kernels for $\mathcal{K}_u$ in the sense that the reproducing property
	$f(\zeta) = \inner{ f, k_{\zeta} }$ holds for all $f$ in $\K_u$.  In a similar manner, the functions
	$$
		[C k_{\zeta}](z) = \frac{u(z) - u(\zeta)}{z - \zeta}
	$$
	are also defined and belong to $\mathcal{K}_u$ whenever $u$ has an ADC at $\zeta$.
	
\subsection{Two Results of Aleksandrov}
	Letting $H^{\infty}$ denote the Banach algebra of all bounded analytic functions on $\D$,
	we observe that the set $\K_u^{\infty} := \K_u \cap H^{\infty}$ is dense in $\K_u$ since
	$\mbox{span}\{ {S^*}^n u : n =1,2,\ldots\}$ is dense in $\K_u$.  Another way to see that $\K_u^{\infty}$
	is dense in $\K_u$ is to observe	
	that each reproducing kernel \eqref{eq-rk-defn} belongs to $\K_u^{\infty}$ whence
	$\mbox{span}\{ k_{\lambda} : \lambda \in \Lambda\}$ is dense in $\K_u$ whenever	
	$\Lambda$ is a uniqueness set for $\K_u$.
	
	For many approximation arguments, the density of $\K_u^{\infty}$ in $\K_u$ is sufficient.
	In other circumstances, however, one requires continuity up to the boundary.  Unfortunately,
	for many inner functions (e.g., singular inner functions)
	it is difficult to exhibit a single nonconstant function in $\K_u$ which is continuous on $\D^{-}$
	(i.e., which belongs to the intersection of $\K_u$ with the \emph{disk algebra} $\A$).  The following
	surprising result asserts that $\K_u \cap \A$, far from being empty, is actually dense in $\K_u$.	
	
	\begin{theorem}[Aleksandrov \cite{MR1359992}]\label{TheoremAleksandrovDensity}
		For $p \in (1,\infty)$, $\K_u^p \cap \A$ is dense in $\K_u^p$.
	\end{theorem}
	
	A detailed exposition of Aleksandrov's density theorem can be found in \cite[p.~186]{CMR}.
	For related results concerning whether or not $\K_u$ contains functions of varying degrees of smoothness,
	the reader is invited to consult \cite{MR2198372}.
	One consequence of Theorem \ref{TheoremAleksandrovDensity} is that it allows us to discuss
	whether or not $\mathcal{K}_{u}^{p}$ can be embedded in $L^{p}(\mu)$ where $\mu$ is a measure on $\dD$.

	\begin{theorem}[Aleksandrov \cite{MR1359992}] \label{A-embedding-T}
		Let $u$ be an inner function, $\mu$ be a positive Borel measure on $\dD$, and $p \in (1, \infty)$.
		If there exists a $C > 0$ such that
		\begin{equation} \label{A-embed}\qquad\qquad\qquad
			\|f\|_{L^p(\mu)} \leqslant C \|f\|_p, \qquad \forall f \in \A \cap \mathcal{K}_{u}^{p},
		\end{equation}
		then every function in $\mathcal{K}_{u}^{p}$ has a finite nontangential limit
		$\mu$-almost everywhere and \eqref{A-embed} holds for all $f \in \mathcal{K}_{u}^{p}$.
	\end{theorem}

	It is clear that every measure on $\dD$
	which is also Carleson measure (see \cite{Garnett})  satisfies \eqref{A-embed}.  However, there are generally
	many other measures which also satisfy \eqref{A-embed}.
	For example, if $u$ has an ADC at $\zeta$, then
	the point mass $\delta_\zeta$ satisfies  \eqref{A-embed} with $p = 2$.
		
\subsection{The Compressed Shift}
	Before introducing truncated Toeplitz operators in general in Subsection \ref{SubsectionTTO},
	we should first introduce and familiarize ourselves with the most important and well-studied example.
	The so-called \emph{compressed shift} operator is simply the compression of the unilateral shift
	\eqref{eq-ShiftOperator} to a model space $\K_u$:
	\begin{equation}\label{eq-CompressedShift}
		A_{z}^{u} := P_{u}S|_{\K_u}.
	\end{equation}
	The adjoint of $A_z^u$ is the restriction of the backward shift \eqref{eq-BackwardShift} to $\K_u$.
	Being the compression of a contraction, it is clear that $A_z^u$ is itself a contraction and in fact,
	such operators and their vector-valued analogues can be used to model certain types of contractive operators
	\cite{Bercovici, N1,N2,N3}.  The following basic properties of $A_z^u$ are well-known and can be found,
	for instance, in \cite{N1,Sarason}.

	\begin{theorem}
		If $u$ is a nonconstant inner function, then
		\begin{enumerate}\addtolength{\itemsep}{0.5\baselineskip}
			\item The invariant subspaces of $A_{z}^{u}$ are $v H^2 \cap (u H^2)^{\perp}$,
				where $v$ is an inner function which divides $u$ (i.e., $u/v \in H^{\infty}$). 	
			\item $A_z^u$ is cyclic with cyclic vector $k_0$.  That is to say, the closed linear span of
			$\{(A^{u}_{z})^{n} k_0: n = 0, 1, 2, \cdots\} = \K_{u}$. Moreover,  $f \in \K_{u}$ is cyclic for
				$A^{u}_{z}$ if and only if $u$ and the inner factor of $f$ are relatively prime.
			\item $A_z^u$ is irreducible (i.e., has no proper, nontrivial reducing subspaces).
		\end{enumerate}
	\end{theorem}

	To discuss the spectral properties of $A_{z}^{u}$ we require the following definition.

	\begin{definition}
		If $u$ is an inner function, then the \emph{spectrum} $\sigma(u)$ of  $u$ is the set
		\begin{equation*}
			\sigma(u) := \left\{\lambda \in \D^{-}: \liminf_{z \to \lambda} |u(z)| = 0\right\}.
		\end{equation*}
	%	If $\lambda \in \dD$, then the preceding limit inferior is to be interpreted in the nontangential limiting sense.
	\end{definition}
	
	If $u = b_{\Lambda} s_{\mu}$, where $b$ is a Blaschke product with zero sequence
	$\Lambda = \{\lambda_n\}$ and $s_{\mu}$ is a singular inner function with corresponding
	singular measure $\mu$, then
	\begin{equation*}
		\sigma(u) =\Lambda^{-} \cup {\mbox{supp}} \mu.
	\end{equation*}
	The following related lemma is well-known result of Moeller and we refer the reader
	to \cite[p.~65]{N1} or \cite[p.~84]{CR} for its proof.
	
	\begin{lemma}\label{LemmaContinuation}
		Each function in $\K_u$ can be analytically continued across $\dD\backslash\sigma(u)$.
	\end{lemma}

	An explicit description of the spectrum of the compressed
	shift $A_z^u$ can be found in Sarason's article \cite[Lem.~2.5]{Sarason}, although portions of it date back to the
	work of Liv\v{s}ic and Moeller \cite[Lec.~III.1]{N1}.	

	\begin{theorem} \label{Comp-shift}
		If $u$ is an inner function, then
		\begin{itemize}\addtolength{\itemsep}{0.5\baselineskip}
			\item[(i)]\quad The spectrum $\sigma(A^{u}_{z})$ of $A_{z}^{u}$ is equal to $\sigma(u)$.
			\item[(ii)]\quad The point spectrum of $\sigma_p(A^{u}_{z})$ of $A_{z}^{u}$ is equal to $\sigma(u) \cap \D$.
			\item[(iii)]\quad The essential spectrum $\sigma_e(A^{u}_{z})$ of $A_{z}^{u}$ is equal to $\sigma(u) \cap \dD$.
		\end{itemize}
	\end{theorem}

\subsection{Clark Unitary Operators and their Spectral Measures}
	Maintaining the notation and conventions of the preceding subsection, let us define,
	for each $\alpha \in \dD$, the following operator on $\K_u$:
	\begin{equation}\label{eq-ClarkOperator}
		U_{\alpha} := A_{z}^u + \frac{\alpha}{1 - \overline{u(0)} \alpha} k_{0} \otimes Ck_0.
	\end{equation}
	In the above, the operator $f \otimes g$, for $f,  g \in H^2$, is given by the formula
	$$(f \otimes g)(h) = \langle h, g\rangle f.$$
	A seminal result of Clark \cite{MR0301534} asserts that each $U_{\alpha}$ is a cyclic unitary operator and, moreover,
	that every unitary, rank-one perturbation of $A^{u}_z$ is of the form \eqref{eq-ClarkOperator}.  Furthermore,
	Clark was even able to concretely identify the corresponding spectral measures $\sigma_{\alpha}$ for these
	so-called \emph{Clark operators}.  We discuss these results below (much of this material is presented
	in greater detail in the recent text \cite{CMR}).

	\begin{theorem}[Clark]
		For each $\alpha \in \dD$, $U_{\alpha}$ is a cyclic unitary operator on $\mathcal{K}_u$.
		Moreover, any unitary rank-one perturbation of $A^{u}_z$ is equal to $U_{\alpha}$ for some $\alpha \in \dD$.
	\end{theorem}

	The spectral theory for the \emph{Clark operators} $U_{\alpha}$ is well-developed and explicit.
	For instance, if $u(0) = 0$, then a point $\zeta$ on $\dD$ is an eigenvalue of $U_{\alpha}$
	if and only if $u$ has an ADC at $\zeta$ and $u(\zeta) = \alpha$. The corresponding eigenvector is
	\begin{equation*}
		k_{\zeta}(z) = \frac{1 - \overline{\alpha} u(z)}{1 - \overline{\zeta} z},
	\end{equation*}
	which is simply a boundary kernel \eqref{eq-BoundaryKernel}.

	Since each $U_{\alpha}$ is cyclic, there exists a measure $\sigma_{\alpha}$, supported on
	$\dD$, so that $U_{\alpha}$ is unitarily equivalent to the operator $M_z:L^2(\sigma_{\alpha})\to L^2(\sigma_{\alpha})$
	of multiplication by the independent variable on the Lebesgue space $L^2(\sigma_{\alpha})$, i.e., $M_z f = z f$.  To concretely
	identify the spectral measure $\sigma_{\alpha}$, we require an
	old theorem of Herglotz, which states that any positive harmonic function  on $\D$ can be written as
	the Poisson integral
	\begin{equation*}
		(\mathfrak{P} \sigma)(z) = \int_{\dD} \frac{1 - |z|^2}{|\zeta - z|^2}\, d\sigma(\zeta)
	\end{equation*}
	of some unique finite, positive, Borel measure $\sigma$ on $\dD$ \cite[Thm.~1.2]{Duren}.

	\begin{theorem}[Clark]
		If $\sigma_{\alpha}$ is the unique measure on $\dD$ satisfying
				\begin{equation} \label{Clark-meas-defn}
			\frac{1 - |u(z)|^2}{|\alpha - u(z)|^2} =  \int_{\dD} \frac{1 - |z|^2}{|\zeta - z|^2} \,d\sigma_{\alpha}(\zeta),
		\end{equation}
		then $U_{\alpha}$ is unitarily equivalent to the operator
		$M_z:L^2(\sigma_{\alpha})\to L^2(\sigma_{\alpha})$ defined by $M_{z} f = z f$.
	\end{theorem}

	The \emph{Clark measures} $\{\sigma_{\alpha}: \alpha \in \dD\}$ corresponding to an inner function
	$u$ have many interesting properties.  We summarize some of these results in the following theorem.
	The reader may consult \cite{CMR} for further details.

	\begin{theorem} \label{Alek-Dec-Form}
		\hfill
		\begin{itemize}\addtolength{\itemsep}{0.5\baselineskip}
			\item[(i)]\quad $\sigma_{\alpha}$ is singular with respect to Lebesgue measure for each $\alpha \in \dD$.
			\item[(ii)]\quad  $\sigma_{\alpha} \perp \sigma_{\beta}$ when $\alpha \not = \beta$.
			\item[(iii)]\quad (Nevanlinna) $\sigma_{\alpha}(\{\zeta\}) > 0$ if and only if
				 $u(\zeta) = \alpha$ and $u$ has an ADC at $\zeta$. Moreover,
				 \begin{equation*}
					\sigma_{\alpha}(\{\zeta\}) = \frac{1}{|u'(\zeta)|}.
				 \end{equation*}
	
			\item[(iv)]\quad (Aleksandrov) For any $f \in C(\dD)$ we have
				\begin{equation} \label{ADF}
					\int_{\dD} \left( \int_{\dD} g(\zeta) \,d\sigma_{\alpha}(\zeta)\right) \!dm(\alpha) = \int_{\dD} g(\zeta) \,dm(\zeta).
				\end{equation}
		\end{itemize}
	\end{theorem}
	
	Condition (iv) of the preceding theorem is a special case of the \emph{Aleksandrov disintegration theorem}:
	If $g$ belongs to $L^1$, then the map
	\begin{equation*}
		\alpha \to \int_{\dD} g(\zeta) \,d\sigma_{\alpha}(\zeta)
	\end{equation*}
	is defined for $m$-almost every $\alpha$ in $\dD$ and, as a function of $\alpha$, it belongs to $L^1$
	and satisfies the natural analogue of \eqref{ADF}.  In fact, the
	Clark measures $\sigma_{\alpha}$ are often called \emph{Aleksandrov-Clark} measures
	in light of Aleksandrov's deep work on the subject, which actually generalizes to measures $\mu_{\alpha}$ on $\dD$
	satisfying
	\begin{equation*}
		\frac{1 - |u(z)|^2}{|\alpha - u(z)|^2} =  \int_{\dD} \frac{1 - |z|^2}{|\zeta - z|^2} \,d\mu_{\alpha}(\zeta)
	\end{equation*}
	for arbitrary functions $u$ belonging to the unit ball of $H^{\infty}$.  The details and ramifications of
	this remarkable result are discussed in detail in \cite{CMR}.

\subsection{Finite Dimensional Model Spaces}
	It is not hard to show that the model space $\K_u$ is finite dimensional if and only if
	$u$ is a finite Blaschke product.  In fact, if $u$ is a finite Blaschke product with zeros
	$\lambda_1, \lambda_2, \ldots, \lambda_N$, repeated according to multiplicity,
	then $\dim \K_u = N$ and
	\begin{equation}\label{eq-PolynomialSpace}
		\K_u = \left\{\frac{\sum_{j = 0}^{N - 1} a_j z^j}{\prod_{j = 1}^{N} (1 - \overline{\lambda_{j}} z)}: a_{j} \in \C\right\}.
	\end{equation}
	With respect to the representation \eqref{eq-PolynomialSpace},
	the conjugation \eqref{eq-ModelConjugation} on $\K_u$ assumes the simple form
	\begin{equation*}
		C\left(\frac{\sum_{j = 0}^{N - 1} a_j z^j}{\prod_{j = 1}^{N} (1 - \overline{\lambda_{j}} z)}\right)
		= \frac{\sum_{j = 0}^{N - 1} \overline{a_{N-1-j}} z^j}{\prod_{j = 1}^{N} (1 - \overline{\lambda_{j}} z)}.
	\end{equation*}
	If the zeros of $u$ are distinct, then the Cauchy kernels $c_{\lambda_i}$ from \eqref{eq-CauchyKernel}
	corresponding to the $\lambda_i$ form a basis for $\K_u$ whence
	\begin{equation*}
		\K_u = \mbox{span}\{ c_{\lambda_1}, c_{\lambda_2}, \ldots, c_{\lambda_N}\}.
	\end{equation*}
	If some of the $\lambda_i$ are repeated, then one must include the appropriate derivatives of the $c_{\lambda_i}$ to obtain
	a basis for $\K_u$.
	
	Although the natural bases for $\K_u$ described above are not orthogonal,
	a particularly convenient orthonormal basis exists.  For $\lambda \in \mathbb{D}$, let
	\begin{equation*}
		b_{\lambda}(z) = \frac{z - \lambda}{1 - \overline{\lambda} z}
	\end{equation*}
	be a disk automorphism with a zero at $\lambda$ and for each $1 \leq n \leq N$ let
	\begin{equation} \label{basis-Tak}
		\gamma_{n}(z) = \frac{\sqrt{1 - |\lambda_n|^2}}{1 - \overline{\lambda_{n}} z} \prod_{k = 1}^{n - 1} b_{\lambda_{k}}(z).
	\end{equation}
	The following important fact was first observed by Takenaka \cite{Tak} in 1925, although it
	has been rediscovered many times since.

	\begin{theorem}[Takenaka]
		$\{\gamma_1,\gamma_2,\ldots,\gamma_N\}$ is an orthonormal basis for $\K_{u}$.
	\end{theorem}

	If $u$ is an infinite Blaschke product, then an extension, due to Walsh \cite{N1},
	of the preceding result tells us that $\{\gamma_1,\gamma_2,\ldots\}$ is an
	orthonormal basis for $\K_{u}$.

	Let us return now to the finite dimensional setting.  Suppose that $u$ is
	a finite Blaschke product with $N$ zeros, repeated according
	to multiplicity.  To avoid some needless technical details, we assume that $u(0) = 0$.
	If $\{\zeta_1,\zeta_2, \ldots, \zeta_N\}$ are the eigenvalues of the Clark unitary operator
	\begin{equation*}
		U_{\alpha} := A^{u}_{z} + \alpha k_{0} \otimes Ck_0,
	\end{equation*}
	(i.e., the $N$ distinct solutions on $\dD$ to $u(\zeta) = \alpha$),
	then the corresponding eigenvectors $\{k_{\zeta_1}, k_{\zeta_2},\ldots, k_{\zeta_N}\}$
	are orthogonal.  A routine computation shows that $\|k_{\zeta_j}\| = \sqrt{|u'(\zeta)|}$ so that
	\begin{equation} \label{basis-Clark}
		\left\{\frac{k_{\zeta_1}}{\sqrt{|u'(\zeta_1)|}},\frac{k_{\zeta_2}}{\sqrt{|u'(\zeta_2)|}},
		\ldots, \frac{k_{\zeta_N}}{\sqrt{|u'(\zeta_N)|}}\right\}
	\end{equation}
	is an orthonormal basis for $\K_u$.  This is called a \emph{Clark basis} for $\K_u$.
	Letting $w_j = \exp\left(-\frac{1}{2}(\arg\zeta_j - \arg\alpha)\right)$, it turns out that
	\begin{equation} \label{basis-modClark}
		\left\{\frac{w_1 k_{\zeta_1}}{\sqrt{|u'(\zeta_1)|}}, \frac{w_2 k_{\zeta_2}}{\sqrt{|u'(\zeta_2)|}},
		\ldots, \frac{w_N k_{\zeta_N}}{\sqrt{|u'(\zeta_N)|}}\right\}
	\end{equation}
	is an orthonormal basis for $\K_u$, each vector of which is fixed by the conjugation
	\eqref{eq-ModelConjugation} on $\K_u$ (i.e., in the terminology of \cite{G-P}, \eqref{basis-modClark}
	is a \emph{$C$-real} basis for $\K_u$).  We refer to a basis of the form \eqref{basis-modClark} as
	a \emph{modified Clark basis} for $\K_u$ (see \cite{CCO} and \cite{NLEPHS} for further details).

%%%%%%%%%%%%%%%%%%%%%%%%%
\subsection{Truncated Toeplitz Operators}\label{SubsectionTTO}
	The \emph{truncated Toeplitz operator} $A_{\phi}^u$ on $\K_u$ having \emph{symbol} $\phi$ in $L^2$
	is the closed, densely defined operator
	\begin{equation*}
		A^{u}_{\phi} f := P_u (\phi f)
	\end{equation*}
	having domain\footnote{Written as an integral transform,
	$P_u$ can be regarded as an operator from $L^1$ into $\mbox{Hol}(\D)$.}
	\begin{equation*}
		\mathscr{D}(A^{u}_{\phi}) = \{f \in \K_{u}: P_u(\phi f) \in \K_u\}.
	\end{equation*}
	When there is no danger of confusion, we sometimes write $A_{\phi}$ in place of $A_{\phi}^u$.
	A detailed discussion of unbounded truncated Toeplitz operators and their properties can be found
	in Section \ref{SectionUnbounded}.  For the moment, we focus on those truncated Toeplitz operators
	which can be extended to bounded operators on $\K_u$.
	
	\begin{definition}
		Let $\TTO_u$ denote the set of all {\em bounded} truncated Toeplitz operators on $\K_u$.
	\end{definition}
	
	For Toeplitz operators, recall that $\norm{T_{\phi}} = \norm{\phi}_{\infty}$ holds for
	each $\phi$ in $L^{\infty}$.  In contrast, we can say little more than
	\begin{equation}\label{eq-TrivialEstimate}
		0 \leqslant \norm{A^{u}_{\phi}} \leqslant \norm{\phi}_{\infty}
	\end{equation}
	for general truncated Toeplitz operators.
	In fact, computing, or at least estimating, the norm of a truncated Toeplitz operator is a difficult problem.
	This topic is discussed in greater detail in Section \ref{SectionNorm}.
	However, a complete characterization of those symbols which yield the zero operator
	has been obtained by Sarason \cite[Thm.~3.1]{Sarason}.	
		
	\begin{theorem}[Sarason]\label{TheoremNonunique}
		A truncated Toeplitz operator $A^{u}_{\phi}$ is identically zero if and only if $\phi \in uH^2 + \overline{u H^2}$.
	\end{theorem}
	
	In particular, the preceding result tells us that there are always infinitely many symbols (many of them unbounded)
	which represent the same truncated Toeplitz operator.  On the other hand, since $A_{\phi}^u = A_{\psi}^u$
	if and only if $\psi = \phi + u H^2 + \overline{u H^2}$, we actually enjoy some freedom in specifying
	the symbol of a truncated Toeplitz operator.  The following corollary makes this point concrete.

	\begin{corollary}
		If $A$ belongs to $\TTO_u$, then there exist $\phi_1$ and $\phi_2$ in $\K_u$ such that
		$A = A_{\phi_1 + \overline{\phi_2}}$.  Furthermore, $\phi_1$ and $\phi_2$ are uniquely determined
		if we fix the value of one of them at the origin.
	\end{corollary}
	
	To some extent, the preceding corollary can be reversed.
	As noted in \cite{BCFMT}, if we assume that $A \in \TTO_u$ has a symbol $\phi_1 + \overline{\phi_2}$, where
	$\phi_1$ and $\phi_2$ belong to $\K_u$ and $\phi_2(0) = 0$, then we can recover $\phi_1$ and $\phi_2$
	by knowing the action of $A$ on the reproducing kernels $k_{\lambda}$ and the conjugate reproducing kernels
	$C k_{\lambda}$.  Indeed, one just needs to solve the following linear system
	in the variables $\phi_1(\lambda)$ and $\overline{\phi_2(\lambda)}$:
	\begin{align*}
		\phi_1(\lambda) - \overline{u(0)} u(\lambda) \overline{\phi_2(\lambda)} &= \inner{ A k_0, k_{\lambda} }, \\
		\overline{\phi_2(\lambda)} - u(0) \overline{u(\lambda)} \phi_1(\lambda) &= \inner{A C k_0, C k_0 } - \inner{A k_0, k_0 }.
	\end{align*}
	With more work, one can even obtain an estimate of $\max\{\|\phi_1\|, \|\phi_2\|\}$ \cite{BCFMT}.

	Letting $C$ denote the conjugation \eqref{eq-ModelConjugation} on $\K_u$, a direct computation
	confirms the following result from \cite{G-P}:

	\begin{theorem}[Garcia-Putinar]\label{TTOareCSO}
		For any $A \in \TTO_u$, we have $A = C A^* C$.
	\end{theorem}

	In particular, Theorem \ref{TTOareCSO} says each truncated Toeplitz operator is a \emph{complex symmetric operator},
	a class of Hilbert space operators which has undergone much recent
	study \cite{Chevrot, CCO, G-P, G-P-II, ESCSO,VPSBF,ATCSO,ONCPCSO,OCCSO,
	SNCSO, CSPI, Gilbreath, JKLL,JKLL2, ZL, ZLZ, Zag,WXH,Wang, ZLJ}.
	In fact, it is suspected that truncated Toeplitz operators might serve as some sort of model operator for various
	classes of complex symmetric operators (see Section \ref{SectionCSO}). For the moment, let us simply note that
	the matrix representation of a truncated Toeplitz operator $A_{\phi}^u$ with respect to a modified
	Clark basis \eqref{basis-modClark} is \emph{complex symmetric} (i.e., self-transpose).  This was first observed
	in \cite{G-P} and developed further in \cite{CCO}.

	An old theorem of Brown and Halmos \cite{BH} says that a bounded operator $T$ on $H^2$
	is a Toeplitz operator if and only if $T = S T S^{*}$.  Sarason recently obtained a version
	of this theorem for truncated Toeplitz operators \cite[Thm.~4.1]{Sarason}.
	
	\begin{theorem}[Sarason]
		A bounded operator $A$ on $\K_u$ belongs to $\TTO_u$ if and only if there are
		functions $\phi, \psi \in \K_u$ such that
		\begin{equation*}
			A = A^{u}_{z} A (A^{u}_{z})^{*} + \phi \otimes k_0 + k_0 \otimes \psi.
		\end{equation*}
	\end{theorem}

	When $\K_u$ is finite dimensional, one can get more specific results using matrix representations.
	For example, if $u = z^N$, then $\{1, z, \ldots, z^{N - 1}\}$ is an orthonormal basis for $\K_{z^N}$.  Any operator
	in $\TTO_{z^N}$ represented with respect to this basis yields a Toeplitz matrix and, conversely, any
	$N \times N$ Toeplitz matrix gives rise to a truncated Toeplitz operator on $\K_{z^N}$.
	Indeed the matrix representation of $A_{\phi}^{z^N}$ with respect to $\{1,z, \ldots z^{N - 1}\}$
	is the Toeplitz matrix $(\widehat{\phi}(j - k))_{j, k=0}^{N-1}$.
	For more general finite Blaschke products we have the following result from \cite{MR2440673}.
	
	\begin{theorem}[Cima-Ross-Wogen]
		Let $u$ be a finite Blaschke product of degree $n$ with distinct zeros
		$\lambda_1, \lambda_2, \ldots, \lambda_{n}$ and let $A$ be any linear transformation on the $n$-dimensional space
		$\K_u$.  If $M_{A} = (r_{i, j})_{i,j=1}^n$ is the matrix representation of $A$ with respect to the basis
		$\{k_{\lambda_1},k_{\lambda_2}, \cdots, k_{\lambda_{n}}\}$, then $A \in \TTO_u$ if and only if
		\begin{equation*}
			r_{i, j} = \overline{\left(\frac{u'(\lambda_1)}{u'(\lambda_i)}\right)} \left(\frac{r_{1, i} \overline{(\lambda_1 - \lambda_i)}
			+ r_{1, j} \overline{(\lambda_j - \lambda_1)}}{\overline{\lambda_j - \lambda_i}}\right),
		\end{equation*}
		for $1 \leqslant i, j \leqslant n$ and $i \neq j$
	\end{theorem}
	
	Although the study of general truncated Toeplitz operators appears to be difficult, there is a distinguished
	subset of these operators which are remarkably tractable.  We say that $A_{\phi}^u$ is an \emph{analytic}
	truncated Toeplitz operator if the symbol $\phi$ belongs to $H^{\infty}$, or more generally, to $H^2$.  	It turns out
	that the natural polynomial functional calculus $p(A_z^u) = A_p^u$ can be extended to $H^{\infty}$ in such a way
	that the symbol map $\phi \mapsto \phi(A_z^u) := A_{\phi}^u$ is linear, contractive, and multiplicative.
	As a broad generalization of Theorem \ref{Comp-shift}, we have the following spectral mapping theorem
	\cite[p.~66]{N1}, the proof of which depends
	crucially on the famous Corona Theorem of L.~Carleson \cite{MR0141789}.
	
	\begin{theorem}\label{TheoremSpectralMapping}
		If $\phi \in H^{\infty}$, then
		\begin{enumerate}\addtolength{\itemsep}{0.5\baselineskip}
			\item $\sigma(A_{\phi}^u) = \left\{ \lambda : \inf_{z\in\D}(|u(z)| + | \phi(z) - \lambda|)=0\right\}$.
			\item If $\phi \in H^{\infty} \cap C(\dD)$, then $\sigma(A_{\phi}^u) = \phi(\sigma(u))$.
		\end{enumerate}	
	\end{theorem}

	We conclude this section by remarking that vector-valued analogues are available
	for most of the preceding theorems.  However, these do not concern us here and
	we refer the reader to \cite{N1} for further details.

\section{$\TTO_u$ as a Linear Space}\label{SectionLinearSpace}

	Recent work of Baranov, Bessonov, and Kapustin \cite{BBK} has shed significant
	light on the structure of $\TTO_u$ as a linear space.  Before describing these results,
	let us first recount a few important observations due to Sarason.  The next theorem is
	\cite[Thm.~4.2]{Sarason}.

	\begin{theorem}[Sarason]\label{TheoremWOT}
		$\TTO_u$ is closed in the weak operator topology.
	\end{theorem}
	
	It is important to note that $\TTO_u$ is not an operator algebra, for the product of truncated Toeplitz operators
	is rarely itself a truncated Toeplitz operator (the precise conditions under which this occurs were found by Sedlock
	\cite{Sed, SedlockThesis}).  On the other hand, $\TTO_u$ contains a number
	of interesting subsets which are algebras.  The details are discussed in Section \ref{SectionAlgebras}, followed
	in Section \ref{SectionCStar} by a brief discussion about $C^*$-algebras generated by
	truncated Toeplitz operators.
		
	In order to better frame the following results, first
	recall that there are no nonzero compact Toeplitz operators on $H^2$ \cite{BH}.
	In contrast, there are many examples of finite rank (hence compact) truncated Toeplitz operators.
	In fact, the rank-one truncated Toeplitz operators were first identified by Sarason \cite[Thm.~5.1]{Sarason}.
	
	\begin{theorem}[Sarason] \label{finite-rank-Sarason}
		For an inner function $u$, the operators
		\begin{enumerate}\addtolength{\itemsep}{0.5\baselineskip}
			\item $k_{\lambda} \otimes C k_{\lambda} = A^u_{\frac{\overline{u}}{\overline{z} - \overline{\lambda}}}$
				for $\lambda \in \D$,
			\item $C k_{\lambda}  \otimes k_{\lambda}= A^u_{\frac{u}{z - \lambda}}$ for $\lambda \in \D$,
			\item $k_{\zeta} \otimes k_{\zeta} = A^u_{k_{\zeta} + \overline{k_{\zeta}} - 1}$
				where $u$ has an ADC at $\zeta \in \dD$,
		\end{enumerate}
		are truncated Toeplitz operators having rank one.
		Moreover, any truncated Toeplitz operator of rank one is a scalar multiple of
		one of the above.
	\end{theorem}

	We should also mention the somewhat more involved results of Sarason \cite[Thms.~6.1 \& 6.2]{Sarason} which identify
	a variety of natural finite rank truncated Toeplitz operators.  Furthermore,  the following linear algebraic description of
	$\TTO_u$ has been obtained \cite[Thm.~7.1]{Sarason} in the finite dimensional setting.

	\begin{theorem}[Sarason]
		If $\dim \K_u = n$, then
		\begin{enumerate}\addtolength{\itemsep}{0.5\baselineskip}
			\item $\dim \TTO_u = 2n-1$,
			\item If $\lambda_1, \lambda_2, \ldots, \lambda_{2n-1}$ are distinct points of $\D$, then the operators
			$k_{\lambda_j}^u \otimes \tilde{k}_{\lambda_j}^u$ for $j = 1,2,\ldots, 2n-1$ form a basis for $\TTO_u$.\footnote{Recall that we are using the notation $\widetilde{f} := C f$ for $f \in \mathcal{K}_u$. }
		\end{enumerate}
	\end{theorem}

	When confronted with a novel linear space, the first questions to arise concern duality.
	Baranov, Bessonov, and Kapustin recently identified the predual of $\TTO_u$
	and discussed the weak-$*$ topology on $\TTO_u$ \cite{BBK}.  Let us briefly summarize
	some of their major results.  First consider the space
	\begin{equation*}
		\mathcal{X}_u := \left\{F =  \sum_{n = 1}^{\infty} f_n \overline{g_n}\,\,:\,\, f_n, g_n \in \mathcal{K}_u,
		\,\,\sum_{n = 1}^{\infty} \|f_n\| \|g_n\| < \infty \right\}
	\end{equation*}
	with norm
	\begin{equation*}
		\norm{F}_{\mathcal{X}_u} := \inf\left\{\sum_{n = 1}^{\infty} \|f_n\| \|g_n\|: F = \sum_{n = 1}^{\infty} f_n \overline{g_n}\right\}.
	\end{equation*}
	It turns out that
	\begin{equation*}
		\mathcal{X}_u \subseteq \overline{u} z H^1 \cap u \overline{z H^1},
	\end{equation*}
	and that each element of $\mathcal{X}_u$ can be written as a linear combination of four elements of the form
	$f \overline{g} $, where $f$ and $g$ belong to $\mathcal{K}_u$.
	 The importance of the space $\mathcal{X}_u$ lies in the following important theorem and its corollaries.

	\begin{theorem}[Baranov-Bessonov-Kapustin \cite{BBK}]
		For any inner function $u$, $\mathcal{X}_{u}^{*}$, the dual space of $\mathcal{X}_u$, is isometrically
		isomorphic to $\TTO_u$ via the dual pairing
		\begin{equation*}
			(F, A) := \sum_{n = 1}^{\infty} \langle A f_n, g_n \rangle, \quad F = \sum_{n = 1}^{\infty} f_n \overline{g_n}, \quad A \in \TTO_u.
		\end{equation*}
		Furthermore, if $\TTO^{c}_{u}$ denotes the compact truncated Toeplitz operators, then $(\TTO_{u}^{c})^{*}$,
		the dual of $\TTO_{u}^{c}$, is isometrically isomorphic to $\mathcal{X}_u$.
	\end{theorem}

	\begin{corollary}[Baranov-Bessonov-Kapustin]
		\hfill
		\begin{enumerate}\addtolength{\itemsep}{0.5\baselineskip}
			\item The weak topology and the weak-$\ast$ topology on $\TTO_u$ are the same.
			\item The norm closed linear span of the rank-one truncated Toeplitz operators is $\TTO_{u}^{c}$.
			\item $\TTO_{u}^{c}$ is weakly dense in $\TTO_u$.
		\end{enumerate}
	\end{corollary}
	
	For a general inner function $u$, we will see below that not every bounded truncated
	Toeplitz operator on $\K_u$ has a bounded symbol (see Section \ref{SectionBoundedSymbol}).
	On the other hand, the following corollary holds in general.
	
	\begin{corollary}[Baranov-Bessonov-Kapustin] \label{wkly-dense}
		The truncated Toeplitz operators with bounded symbols are weakly dense in $\TTO_u$.
	\end{corollary}

	This leaves open the following question.

	\begin{Question}
		Are the truncated Toeplitz operators with bounded symbols \emph{norm} dense in $\TTO_u$?
	\end{Question}

\section{Norms of Truncated Toeplitz Operators}\label{SectionNorm}
	Recall that for $\phi$ in $L^{\infty}$ we have the trivial estimates \eqref{eq-TrivialEstimate} on the norm of a truncated
	Toeplitz operator, but little other general information concerning this quantity.
	For $\phi$ in $L^2$, we may also consider the potentially unbounded truncated Toeplitz
	operator $A_{\phi}^u$.  Of interest is the quantity
	\begin{equation}\label{eq-Norm}
		\|A_{\phi}^{u}\| := \sup\{\|A^{u}_{\phi} f\|: f \in \mathcal{K}_{u} \cap H^{\infty}, \|f\| = 1\},
	\end{equation}
	which we regard as being infinite if $A_{\phi}^u$ is unbounded.  For $A_{\phi}^u$ bounded, \eqref{eq-Norm}
	is simply the operator norm of $A_{\phi}^u$ in light of Theorem \ref{TheoremAleksandrovDensity}.
	Evaluation or estimation of \eqref{eq-Norm} is further complicated by the fact that
	the representing symbol $\phi$ for $A^{u}_{\phi}$ is never unique (Theorem
	\ref{TheoremNonunique}).

	If $u$ is a finite Blaschke product (so that the corresponding model space $\K_u$ is finite dimensional) and $\phi$
	belongs to $H^{\infty}$, then straightforward residue computations allow us to
	represent $A^{u}_{\phi}$ with respect to any of the
	orthonormal bases mentioned earlier (i.e., the Takenaka \eqref{basis-Tak}, Clark \eqref{basis-Clark},
	or modified Clark  \eqref{basis-modClark} bases).  For $\mathcal{K}_{z^n}$, the Takenaka
	basis is simply the monomial basis $\{1, z, z^2, \ldots, z^{n - 1}\}$ and the matrix representation
	of $A_{\phi}^u$ is just a lower
	triangular Toeplitz matrix.  In any case, one can readily compute the norm of $A^{u}_{\phi}$
	by computing the norm of one of its matrix representations.  This approach was undertaken
	by the authors in \cite{NLEPHS}.  One can also approach this problem using the theory of Hankel operators
	(see \cite[eq.~2.9]{Peller} and the method developed in \cite{MR2597679}).
	
	Let us illustrate the general approach with a simple example.
	If $\phi$ belongs to $H^{\infty}$ and $u$ is the finite Blaschke product
	with distinct zeros $\lambda_1, \lambda_2,\ldots, \lambda_n$, then
	the matrix representation for $A^{u}_{\phi}$ with respect to the
	modified Clark basis \eqref{basis-modClark} is
	\begin{equation} \label{mod-Clark-ex}
		\left( \frac{w_k}{\sqrt{|u'(\zeta_k)|}} \frac{w_j}{\sqrt{|u'(\zeta_j)|}}
		\sum_{i = 1}^{n} \frac{\phi(\lambda_i)}{u'(\lambda_i)
		(1 - \overline{\zeta_k} \lambda_i)(1 - \overline{\zeta_j} \lambda_i)}\right)_{j, k=1}^n.
	\end{equation}
	In particular, observe that this matrix is complex symmetric, as predicted by
	Theorem \ref{TTOareCSO}.  As a specific example, consider the Blaschke product
	\begin{equation*}
		u(z) = z \frac{z - \frac{1}{2}}{1 - \frac{1}{2} z}
	\end{equation*}
	and the $H^{\infty}$ function
	\begin{equation*}
		\phi(z) = \frac{2 z - \frac{1}{2}}{1 - \frac{1}{2} z}.
	\end{equation*}
	The parameters in \eqref{mod-Clark-ex} are
	\begin{equation*}
		\alpha=1,\qquad \lambda_{1} = 0, \,\,\lambda_2 = \tfrac{1}{2},
		\qquad \zeta_1 = 1, \,\,\zeta_2 = -1, \qquad w_1 = 1,\,\, w_2 = -i,
	\end{equation*}
	which yields
	\begin{equation*}
		\|A^{u}_{\phi}\|
		= \left\|
			\begin{pmatrix}
			 \frac{5}{4} & -\frac{7 i}{4 \sqrt{3}}\vspace{3pt} \\
			 -\frac{7 i}{4 \sqrt{3}} & -\frac{13}{12}
			\end{pmatrix}
		\right\|
		= \frac{1}{6}(7 + \sqrt{37})
		\approx 2.1805.
	\end{equation*}

	On a somewhat different note, it is possible to obtain lower estimates of $\|A_{\phi}^u\|$ for general $\phi$ in $L^2$.
	This can be helpful, for instance, in determining whether a given truncated Toeplitz operator is unbounded.
	Although a variety of lower bounds on $\norm{A_{\phi}^u}$ are provided in \cite{GRNorm}, we focus here on perhaps
	the most useful of these.  We first require the \emph{Poisson integral}
	\begin{equation*}
		(\mathfrak{P} \phi)(z) := \int_{\dD} \frac{1 - |z|^2}{|\zeta - z|^2} \phi(\zeta) \,dm(\zeta)
	\end{equation*}
	of a function $\phi$ in $L^1$.  In particular, recall that $\lim_{r \to 1^{-}}(\mathfrak{P}\phi)(r \zeta) = \phi(\zeta)$
	whenever $\phi$ is continuous at a point
	$\zeta \in \dD$ \cite[p.~32]{Hoffman} or more generally, when $\zeta$ is a Lebesgue point of $\phi$.

	\begin{theorem}[Garcia-Ross]
		If $\phi \in L^2$ and $u$ is inner, then
		\begin{equation*}
			\|A^{u}_{\phi}\| \geq \sup\{ |(\mathfrak{P} \phi)(\lambda)|: \lambda \in \D: u(\lambda) = 0\},
		\end{equation*}
		where the supremum is regarded as $0$ if $u$ never vanishes on $\D$.
	\end{theorem}
	
	\begin{corollary}
		If $\phi$ belongs to $C(\dD)$ and $u$ is an inner function whose zeros accumulate
		almost everywhere on $\dD$, then $\norm{A_{\phi}^u} = \norm{\phi}_{\infty}$.
	\end{corollary}
	
	A related result on \emph{norm attaining symbols} can be found in \cite{NLEPHS}.

	\begin{theorem}[Garcia-Ross]
		If $u$ is inner, $\phi\in H^{\infty}$, and $A^{u}_{\phi}$ is compact,
		then $\|A^{u}_{\phi}\| = \|\phi\|_{\infty}$ if and only if $\phi$ is a scalar multiple
		 of the inner factor of a function from $\mathcal{K}_{u}$.
	\end{theorem}	

	It turns out that the norm of a truncated Toeplitz operator can be related to certain
	classical extremal problems from function theory. For the following discussion we require a few general
	facts about complex symmetric operators \cite{G-P, G-P-II, CCO}. Recall that a \emph{conjugation} on a complex
	Hilbert space $\mathcal{H}$ is a map $C: \mathcal{H} \to \mathcal{H}$ which is conjugate-linear,
	involutive (i.e., $C^2 = I$), and isometric (i.e., $\langle C x, C y\rangle = \langle y, x \rangle$).
	A bounded linear operator $T: \mathcal{H} \to \mathcal{H}$ is called \emph{$C$-symmetric} if $T = C T^* C$ and
	\emph{complex symmetric} if there exists a conjugation $C$ with respect to which $T$ is $C$-symmetric.
	Theorem \ref{TTOareCSO} asserts that each operator in $\TTO_u$ is $C$-symmetric with respect to the
	conjugation $C$ on $\mathcal{K}_u$ defined by \eqref{eq-ModelConjugation}.  The following general result from
	\cite{NLEPHS} relates the norm of a given $C$-symmetric operator to a certain extremal problem (as is customary,
	$|T|$ denotes the positive operator $\sqrt{T^*T}$).

	\begin{theorem}[Garcia-Ross] \label{TheoremVariationalNorm}
		If $T:\h\to\h$ is a bounded $C$-symmetric operator, then\hfill
		 \begin{enumerate}\addtolength{\itemsep}{0.5\baselineskip}
			\item $\norm{T} = \sup_{ \norm{x} = 1} | \inner{Tx,Cx} |$.

			\item If $\norm{x} = 1$, then $\|T\| = |\inner{Tx, Cx}|$ if and
				only if $T x = \omega \|T\| C x$ for some unimodular constant $\omega$.
				
			\item If $T$ is compact, then the equation $T x = \|T\| C x$ has
				a unit vector solution. Furthermore, this unit vector solution is unique,
				up to a sign,  if and only if the kernel of the operator $|T| - \|T\| I$ is one-dimensional.
		\end{enumerate}
	\end{theorem}
	
	Applying the Theorem \ref{TheoremVariationalNorm} to $A_{\phi}^u$ we obtain the following result.

	\begin{corollary}
		For inner $u$ and $\phi \in L^{\infty}$
		\begin{equation*}
			\|A^{u}_{\phi}\| = \sup\left\{\left|\frac{1}{2 \pi i} \oint_{\dD} \frac{\phi f^2}{u} d z\right|: f \in \mathcal{K}_u, \|f\| = 1\right\}.
		\end{equation*}
	\end{corollary}

	For $\phi$ in $H^{\infty}$, the preceding supremum can be taken over $H^2$.

	\begin{corollary}
		For inner $u$ and $\phi \in H^{\infty}$
		\begin{equation*}
			\|A^{u}_{\phi}\| = \sup\left\{\left|\frac{1}{2 \pi i} \oint_{\dD} \frac{\phi f^2}{u} d z\right|: f \in H^2, \|f\| = 1\right\}.
		\end{equation*}
	\end{corollary}

	The preceding corollary relates the norm of a truncated Toeplitz operator to a certain \emph{quadratic} extremal problem on $H^2$.
	We can relate this to a classical \emph{linear} extremal problem in the following way.  For a rational function $\psi$ with
	no poles on $\dD$ we have the well-studied classical $H^1$ extremal problem \cite{Duren, Garnett}:
	\begin{equation}\label{eq-LinearProblem}
		\Lambda(\psi) := \sup\left\{\left|\frac{1}{2 \pi i} \oint_{\dD} \psi F\, dz\right|: F \in H^1, \|F\|_{1} = 1\right\}.
	\end{equation}
	On the other hand, basic functional analysis tells us that
	\begin{equation*}
		\Lambda(\psi) = \mbox{dist}(\psi, H^{\infty}).
	\end{equation*}

	Following \cite{NLEPHS}, we recall that the extremal problem $\Lambda(\psi)$ has an \emph{extremal function}
	$F_e$ (not necessarily unique).
	It is also known that $F_e$ can be taken to be outer and hence $F_e  = f^2$ for some $f$ in $H^2$.  Therefore
	the linear extremal problem $\Lambda(\psi)$ and the quadratic extremal problem
	\begin{equation}\label{eq-QuadraticProblem}
		\Gamma(\psi) :=  \sup\left\{\left|\frac{1}{2 \pi i} \oint_{\dD} \psi f^2 dz\right|: f \in H^2, \|f\|= 1\right\}
	\end{equation}
	have the same value.
	The following result from \cite{NLEPHS}, combined with the numerical recipes discussed at the beginning
	of this section, permit one to explicitly evaluate many specific extremal problems.  Before doing so,
	we remark that many of these
	problems can be attacked using the theory of Hankel operators, although in that case one must compute the norm
	of a finite-rank Hankel operator acting on an infinite-dimensional space.  In contrast, the truncated Toeplitz approach
	employs only $n \times n$ matrices.

	\begin{corollary} \label{cor-main-result}
		Suppose that $\psi$ is a rational function
		having no poles on $\dD$ and poles $\lambda_1,\lambda_2, \ldots, \lambda_n$
		lying in $\D$, counted according to multiplicity.
		Let $u$ denote the associated Blaschke product whose zeros are precisely $\lambda_1,\lambda_2, \ldots, \lambda_n$ and
		note that $\phi = u \psi$ belongs to $H^{\infty}$. We then have the following:
		 \begin{enumerate}\addtolength{\itemsep}{0.5\baselineskip}
			\item $\|A^{u}_{\phi}\| = \Gamma(\psi) = \Lambda(\psi)$.

			\item There is a unit vector $f \in \mathcal{K}_{u}$ satisfying
					$A^u_{\phi} f = \|A^u_{\phi}\| C f$
				and any such $f$ is an extremal function for $\Gamma(\psi)$.
				In other words,				\begin{equation*}
					 \left|\frac{1}{2 \pi i} \oint_{\dD}  \psi  f^2  \, dz\right|
					= \|A^u_{\phi}\|.
				\end{equation*}
				
			\item Every extremal function $f$ for $\Gamma(\psi)$ belongs to $\mathcal{K}_u$ and satisfies
				\begin{equation*}
					A^u_{\phi} f = \|A_{ \phi}\| C f.
				\end{equation*}

			\item An extremal function for $\Gamma(\psi)$ is unique, up
				to a sign, if and only if the kernel
				of the operator $|A^u_{\phi}| - \|A^u_{\phi}\| I$ is one-dimensional.
		\end{enumerate}
	\end{corollary}
	
	We refer the reader to \cite{NLEPHS} for several worked examples of classical extremal problems $\Lambda(\psi)$
	along with a computation of several extremal functions $F_e$.   For rational functions $\psi$ with no poles on
	$\dD$, we have seen that the linear \eqref{eq-LinearProblem} and the quadratic \eqref{eq-QuadraticProblem}
	extremal problems have the same value.
	Recent work of Chalendar, Fricain, and Timotin shows that this holds in much greater generality.

	\begin{theorem}[Chalendar-Fricain-Timotin \cite{MR2597679}]\label{TheoremCFTExtremal}
		For each $\psi$ in $L^{\infty}$, $\Gamma(\psi) = \Lambda(\psi)$.
	\end{theorem}

	It is important to note that for general $\psi$ in $L^{\infty}$ an extremal function for $\Lambda(\psi)$ need not exist
	(see \cite{NLEPHS} for a relevant discussion).
	Nevertheless, for $\psi$ in $L^{\infty}$, Chalendar, Fricain, and Timotin prove that $\Lambda(\psi) = \Gamma(\psi)$
	by using the fact that $\Lambda(\psi) = \|H_{\psi}\|$, where $H_{\psi}: H^2 \to L^2 \ominus H^2$ is the corresponding
	\emph{Hankel operator} $H_{\psi} f = P_{-}(\psi f)$.  Here $P_{-}$ denotes the orthogonal projection from
	$L^2$ onto $L^2 \ominus H^2$. We certainly have the inequality
	\begin{equation*}
		\Gamma(\psi) = \Lambda(\psi) = \|H_{\psi}\| = \mbox{dist}(\psi, H^{\infty}) \leqslant \|\psi\|_{\infty}.
	\end{equation*}
	When equality holds in the preceding, we say that the symbol $\psi$ is \emph{norm attaining}.
	The authors of \cite{MR2597679} prove the following.

	\begin{theorem}[Chalendar-Fricain-Timotin]
		If $\psi \in L^{\infty}$ is norm attaining then $\psi$ has constant modulus and there exists an extremal outer function for $\Lambda(\psi)$.
	\end{theorem}
	
	Before proceeding, we should also mention the fact that computing the norm of certain truncated
	Toeplitz operators and solving their related extremal problems have been examined for quite some
	time in the study of $H^{\infty}$ control theory and skew-Toeplitz operators
	\cite{MR945002, MR1639646, MR901235, MR971878}.  In the scalar setting, a \emph{skew-Toeplitz operator}
	is a truncated Toeplitz operator $A^{u}_{\phi}$, where the symbol takes the form
	\begin{equation*}
		\qquad\qquad\phi(\zeta) = \sum_{j, k = 0}^{n} a_{j, k} \zeta^{j - k}, \qquad a_{j, k} \in \mathbb{R},
	\end{equation*}
	making $A_{\phi}^{u}$ self-adjoint. In $H^{\infty}$ control theory, the extremal problem
	$$\mbox{dist}(\psi, u H^{\infty}),$$ where $\psi$ is a rational function belonging to $H^{\infty}$,
	plays an important role. From the preceding results, we observe that
	$\|A^{u}_{\psi}\| = \mbox{dist}(\psi, u H^{\infty})$.

\section{The Bounded Symbol and Related Problems}\label{SectionBoundedSymbol}
	Recall that $\TTO_u$ denotes the set of all truncated Toeplitz operators $A_{\phi}^{u}$, densely defined on
	$\K_u$ and having symbols $\phi$ in $L^2$, that can be extended to bounded operators on all of $\K_u$.
	As a trivial example, if $\phi$ belongs to $L^{\infty}$, then clearly $A_{\phi}^u$ belongs to $\TTO_u$.  A major
	open question involved the converse. In other words, if $A_{\phi}^u$ is a
	bounded truncated Toeplitz operator, does there exist a symbol $\phi_0$ in $L^{\infty}$ such that
	$A_{\phi}^u = A_{\phi_0}^u$?
	This question was recently resolved in the negative by
	Baranov, Chalendar, Fricain, Mashreghi, and Timotin \cite{BCFMT}.
	We describe this groundbreaking work, along with important related contributions by
	Baranov, Bessonov, and Kapustin \cite{BBK}, below.

	For symbols $\phi$ in $H^2$, a complete and elegant answer to the bounded symbol problem
	is available.  In the following theorem, the implication (i) $\Leftrightarrow$ (ii) below
	is due to Sarason \cite{Sarason-NF}.  Condition (iii) is often referred to as the \emph{reproducing kernel thesis}
	for $A_{\phi}^u$.

	\begin{theorem}[Baranov, Chalendar, Fricain, Mashreghi, Timotin \cite{BCFMT}]\label{TheoremBCFMT1}
		For $\phi \in H^2$, the following are equivalent.
		\begin{enumerate}\addtolength{\itemsep}{0.5\baselineskip}
			\item $A_{\phi}^{u} \in \TTO_u$.
			\item $A_{\phi}^{u} = A_{\phi_0}^{u}$ for some $\phi_0 \in H^{\infty}$.
			\item $\sup_{\lambda \in \D} \left\|A_{\phi}^{u} \frac{k_{\lambda}}{\|k_{\lambda}\|}\right\| < \infty$.
		\end{enumerate}
		Furthermore, there exists a universal constant $C > 0$ so that any $A_{\phi}^{u} \in \TTO_u$ with
		$\phi \in H^2$, has a bounded symbol $\phi_0$ such that
		\begin{equation*}
			\|\phi_0\|_{\infty} \leqslant C \sup_{\lambda \in \D} \left\|A_{\phi}^{u} \frac{k_{\lambda}}{\|k_{\lambda}\|}\right\|.
		\end{equation*}
	\end{theorem}

	The following result demonstrates the
	existence of bounded truncated Toeplitz operators with no bounded symbol (a small amount
	of function theory, discussed below, is required to exhibit concrete examples).
	Recall from Theorem \ref{finite-rank-Sarason} that for $\zeta$ in $\dD$, the rank one operator
	$k_{\zeta} \otimes k_{\zeta}$ belongs to $\TTO_u$ if and only if $u$ has an ADC at $\zeta$. Using two technical lemmas from \cite{BCFMT} (Lemmas 5.1 and 5.2), they prove the following theorem.

	\begin{theorem}[Baranov, Chalendar, Fricain, Mashreghi, Timotin] \label{no-bd-symbol}
		If $u$ has an ADC at $\zeta \in \dD$ and $p \in (2, \infty)$, then the following are equivalent:
		\begin{enumerate}\addtolength{\itemsep}{0.5\baselineskip}
			\item $k_{\zeta} \otimes k_{\zeta}$ has a symbol in $L^p$.
			\item $k_{\zeta} \in L^p$.
		\end{enumerate}
		Consequently, if $k_{\zeta} \not \in L^p$ for some $p \in (2, \infty)$, then
		$k_{\zeta} \otimes k_{\zeta}$ belongs to $\TTO_u$ and has no bounded symbol.
	\end{theorem}

	From Theorem \ref{AC-paper} we know that if $u = b_{\Lambda} s_{\mu}$, where $b$ is a Blaschke product
	with zeros $\Lambda = \{\lambda_n\}_{n = 1}^{\infty}$ (repeated according to multiplicity) and $s_{\mu}$ is
	a singular inner function with corresponding singular measure $\mu$, then
	\begin{equation} \label{AC-ex-cond1}
		k_{\zeta} \in H^2
		\iff \sum_{n = 1}^{\infty} \frac{1 - |\lambda|^2}{|\zeta - \lambda_n|^2} + \int \frac{d \mu(\xi)}{|\xi - \zeta|^2} < \infty.
	\end{equation}
	This was extended \cite{AC70, Cohn} to $p \in (1, \infty)$ as follows:
	\begin{equation} \label{AC-ex-cond2}
		k_{\zeta} \in H^p \iff
		\sum_{n = 1}^{\infty} \frac{1 - |\lambda|^2}{|\zeta - \lambda_n|^p} + \int \frac{d \mu(\xi)}{|\xi - \zeta|^p} < \infty.
	\end{equation}
	Based upon this, one can arrange it so that $k_{\zeta}$ belongs to $H^2$ but not to $H^p$ for any $p > 2$.
	This yields the desired example of a bounded truncated Toeplitz operator which cannot be represented
	using a bounded symbol.
	We refer the reader to the article \cite{BCFMT} where the details and further examples are considered.
	
	Having seen that there exist bounded truncated Toeplitz operators which lack bounded symbols,
	it is natural to ask if there exist inner functions $u$ so that \emph{every} operator in
	$\TTO_u$ has a bounded symbol?  Obviously this holds when $u$ is a finite Blaschke product.  Indeed,
	in this case the symbol can be taken to be a polynomial in $z$ and $\overline{z}$.  A more difficult result
	is the following (note that the initial symbol $\phi$ belongs to $L^2$, as opposed to $H^2$, as was the
	case in Theorem \ref{TheoremBCFMT1}).

	\begin{theorem}[Baranov, Chalendar, Fricain, Mashreghi, Timotin]\label{TheoremBCFMT2}
		If $a > 0$, $\zeta \in \dD$, and
		\begin{equation*}
			u(z) = \exp\left(a \frac{z + \zeta}{z - \zeta}\right),
		\end{equation*}
		then the following are equivalent for $\phi \in L^2$:
		\begin{enumerate}\addtolength{\itemsep}{0.5\baselineskip}
			\item $A_{\phi}^{u} \in \TTO_u$.
			\item $A_{\phi}^u = A_{\phi_0}^{u}$ for some $\phi_0 \in L^{\infty}$.
			\item $\sup_{\lambda \in \D} \left\|A_{\phi}^{u} \frac{k_{\lambda}}{\|k_{\lambda}\|}\right\| < \infty$.
		\end{enumerate}
		Furthermore, there exists a universal constant $C > 0$ so that any $A_{\phi}^{u} \in \TTO_u$
		with $\phi \in L^2$, has a bounded symbol $\phi_0$ such that
		\begin{equation*}
			\|\phi_0\|_{\infty} \leqslant C \sup_{\lambda \in \D} \left\|A_{\phi}^{u} \frac{k_{\lambda}}{\|k_{\lambda}\|}\right\|.
		\end{equation*}
	\end{theorem}

	In light of Theorems \ref{TheoremBCFMT1} and \ref{TheoremBCFMT2}, one might wonder whether
	condition (iii) (the reproducing kernel thesis) is always equivalent to asserting that $A_{\phi}^u$
	belongs to $\TTO_u$.  Unfortunately, the answer is again negative \cite[Sec.~5]{BCFMT}.
		
	On a positive note, Baranov, Bessonov, and Kapustin recently discovered
	a condition on the inner function $u$ which ensures that \emph{every} operator in
	$\TTO_u$ has a bounded symbol \cite{BBK}.  After a few preliminary details, we discuss
	their work below.

	\begin{definition}
		For $p > 0$, let $\mathcal{C}_{p}(u)$ denote the finite complex Borel measures
		$\mu$ on $\dD$ such that $\K_{u}^{p}$ embeds continuously into $L^p(|\mu|)$.
	\end{definition}

	Since $S^{*} u$ belongs to $\mathcal{K}_u$, it follows from Aleksandrov's embedding theorem
	(Theorem \ref{A-embedding-T}) that for each $\mu$ in $\mathcal{C}_{2}(u)$, the boundary values
	of $u$ are defined $|\mu|$-almost everywhere. Moreover, it turns out that
	$|u| = 1$ holds $|\mu|$-almost everywhere \cite{MR1359992, BBK}.
	For $\mu \in \mathcal{C}_{2}(u)$ the quadratic form
	\begin{equation*}
		(f, g) \mapsto \int_{\dD} f \overline{g} \,d \mu, \quad f, g \in \mathcal{K}_u
	\end{equation*}
	is continuous and so, by elementary functional analysis, there is a bounded
	operator $\mathcal{A}_{\mu}: \mathcal{K}_{u} \to \mathcal{K}_{u}$ such that
	\begin{equation*}
		\langle \mathcal{A}_{\mu} f, g \rangle  = \int_{\dD} f \overline{g}\, d \mu.
	\end{equation*}
	The following important result of Sarason \cite[Thm.~9.1]{Sarason} asserts that
	each such $\mathcal{A}_{\mu}$ is a truncated Toeplitz operator.
	
	\begin{theorem}[Sarason]
		$\mathcal{A}_{\mu} \in \TTO_u$ whenever $\mu \in \mathcal{C}_{2}(u)$.
	\end{theorem}
	
	A natural question, posed by Sarason \cite[p.~513]{Sarason}, is whether the converse holds.
	In other words, does every bounded truncated Toeplitz operator arise from a so-called
	\emph{$u$-compatible measure} \cite[Sect.~9]{Sarason}?  This question was recently settled in the
	affirmative by Baranov, Bessonov, and Kapustin \cite{BBK}.
	
	\begin{theorem}[Baranov-Bessonov-Kapustin \cite{BBK}]
		$A \in \TTO_u$ if and only if $A = \mathcal{A}_{\mu}$ for some $\mu \in \mathcal{C}_{2}(u)$.
	\end{theorem}

	The measure $\mu$ above is called the \emph{quasi-symbol} for the truncated Toeplitz operator.
	For $\phi \in L^{\infty}$ we adopt the convention that $\mathcal{A}_{\phi\, d m}^u := A_{\phi}^{u}$
	so that every bounded symbol is automatically a quasi-symbol.

	It turns out that $\mathcal{C}_{1}(u^2) \subseteq \mathcal{C}_{2}(u) = \mathcal{C}_2(u^2)$ always holds.
	Baranov, Bessonov, and Kapustin showed that equality is the precise condition which ensures that
	every $A \in \TTO_u$ can be represented using a bounded symbol.
	
	\begin{theorem}[Baranov-Bessonov-Kapustin]
		An operator $A \in \TTO_u$ has a bounded symbol if and only if $A = \mathcal{A}_{\mu}$ for some
		 $\mu \in \mathcal{C}_{1}(u^2)$.  Consequently, every operator in $\TTO_u$ has a bounded
		 symbol if and only if $\mathcal{C}_{1}(u^2) = \mathcal{C}_{2}(u)$.
	\end{theorem}

	Recall that each function $F$ in $H^1$ can be written as the product $F = fg$ of two functions in $H^2$.
	Conversely, the product of any pair of functions in $H^2$ lies in $H^1$.  For each $f, g \in \mathcal{K}_u$ we note that
	\begin{equation*}
		H^1 \ni fg = \tilde{\tilde{f}}\tilde{\tilde{g}} = \overline{ \tilde{f} \tilde{g}z^2} u^2 \in \overline{z} u^2 \overline{H_0^1},
	\end{equation*}
	whence $f g \in H^1 \cap \overline{z} u^2 \overline{H^{1}_{0}}$.
	Moreover, one can show that finite linear combinations of pairs of products of functions from $\K_u$
	form a dense subset of  $H^1 \cap \overline{z} u^2 \overline{H^{1}_{0}}$. As a consequence,
	this relationship between $\K_u$ and $H^1 \cap \overline{z} u^2 \overline{H^{1}_{0}}$ is sometimes denoted
	\begin{equation*}
		\mathcal{K}_{u} \odot \mathcal{K}_{u} = H^1 \cap \overline{z} u^2 \overline{H^{1}_{0}}.
	\end{equation*}
	For certain inner functions, one can say much more.

	\begin{theorem}[Baranov-Bessonov-Kapustin]
		For an inner function $u$ the following statements are equivalent.
		\begin{enumerate}\addtolength{\itemsep}{0.5\baselineskip}
			\item $\mathcal{C}_{1}(u^2) = \mathcal{C}_2(u)$.
			\item For each $f \in H^1 \cap \overline{z} u^2 \overline{H^{1}_{0}}$ there exists
				sequences
				$g_j, h_j$ in $\mathcal{K}_u$ such that $\sum_{j} \|g_j\| \|h_j\| < \infty$ and
				\begin{equation*}
					f = \sum_{j} g_j h_j.
				\end{equation*}
				Moreover, there exists a universal $C > 0$, independent of $f$, such that the $g_j, h_j$ can be chosen to satisfy
				$\sum_{j} \|g_j\| \|h_j\| \leqslant C \|f\|_1$.
						\end{enumerate}
	\end{theorem}

	As we have seen, the condition $\mathcal{C}_{1}(u^2)  = \mathcal{C}_{2}(u)$ is of primary importance.
	Unfortunately, it appears difficult to test whether or not a given inner function has this property.  On the other hand,
	the following related geometric condition appears somewhat more tractable.

	\begin{definition}
		An inner function $u$ is called a \emph{one-component inner function} if the set
		$\{z \in \D: |u(z)| < \epsilon\}$ is connected for some $\epsilon > 0$.
	\end{definition}

	One can show, for instance, that for the atomic inner functions $s_{\delta_{\zeta}}$ (where $\delta_{\zeta}$ denotes
	the point mass at a point $\zeta$ on $\dD$), the set $\{|u| < \epsilon\}$ is a disk internally tangent to $\dD$ at $\zeta$.
	In other words, these inner functions
	are one-component inner functions.  The relevance of one-component inner functions
	lies in the following result of Aleksandrov.

	\begin{theorem}[Aleksandrov \cite{MR1734326}]
		If $u$ is inner, then
		\begin{enumerate}\addtolength{\itemsep}{0.5\baselineskip}
			\item $u$ is a one-component inner function if and only if
				\begin{equation*}
					\sup_{\lambda \in \D} \frac{\|k_{\lambda}\|_{\infty}}{\|k_{\lambda}\|} < \infty.
				\end{equation*}
			\item If $u$ is a one-component inner function, then
				 $\mathcal{C}_{p_1}(u) = \mathcal{C}_{p_2}(u)$ for all $p_1, p_2 > 0$.
		\end{enumerate}
	\end{theorem}
	
	There is also a related result of Treil and Volberg \cite{TV}. We take a moment to mention that W. Cohn \cite{Cohn} has described the set $\mathcal{C}_{2}(u)$ in the case where
	$u$ is a one component inner function.
	
	\begin{theorem}[Cohn]
		If $u$ is a one component inner function and $\mu$ is a positive measure on $\dD$
		for which $\mu(\sigma(u) \cap \dD) = 0$, then $\mu \in \mathcal{C}_{2}(u)$
		if and only if there is a constant $C > 0$ such that
		\begin{equation*}
			\int_{\dD} \frac{1 - |z|^2}{|\zeta - z|^2}\,d \mu(\zeta) \leqslant \frac{C}{1 - |u(z)|^2}
		\end{equation*}
		for all $z$ in $\D$.
	\end{theorem}
	
	Note that if $u$ is a one-component inner function, then so is $u^2$.
	Combining the preceding results
	we obtain the following.

	\begin{corollary}
		If $u$ is a one-component inner function, then every operator in $\TTO_u$ has a bounded symbol.
	\end{corollary}
	
	Of course the natural question now (and conjectured in \cite{BCFMT}) is the following:

	\begin{Question}
		Is the converse of the preceding corollary true?
	\end{Question}

	It turns out that there is an interesting and fruitful interplay between the material discussed above and the family
	of Clark measures $\{\sigma_{\alpha}: \alpha \in \dD\}$, defined by \eqref{Clark-meas-defn},
	associated with an inner function $u$.  To be more specific, Aleksandrov showed
	in \cite{MR1359992} that if $\mathcal{C}_{1}(u^2) = \mathcal{C}_{2}(u^2)$ (the equivalent
	condition for every operator in $\TTO_u$ to have a bounded symbol), then every Clark measure is discrete.
	This leads to the following corollary.

	\begin{corollary}
		Let $u$ be an inner function.
		If for some $\alpha \in \dD$ the Clark measure $\sigma_{\alpha}$ is not discrete,
		then there is an operator in $\TTO_u$ without a bounded symbol.
	\end{corollary}

	Since \emph{any} singular measure $\mu$ (discrete or not) is equal to $\sigma_1$
	for some inner function $u$ \cite{CMR}, it follows that if we let $\mu$ be a continuous singular measure, then
	the corollary above yields an example of a truncated Toeplitz operator space $\TTO_u$ which contains
	operators without a bounded symbol.

\section{The Spatial Isomorphism Problem}\label{SectionSpatialIsomorphism}
	For two inner functions $u_1$ and $u_2$, when is $\TTO_{u_1}$ \emph{spatially isomorphic} to $\TTO_{u_2}$?
	In other words, when does there exist a unitary operator $U: \K_{u_1} \to \K_{u_2}$ such that
	$U \TTO_{u_1} U^{*}  = \TTO_{u_2}$?  This is evidently a stronger condition than isometric isomorphism since one insists
	that the isometric isomorphism is implemented in a particularly restrictive manner.

	A concrete solution to the spatial isomorphism problem posed above was given in \cite{TTOSIUES}.  Before
	discussing the solution, let us briefly introduce three basic families of
	spatial isomorphisms between truncated Toeplitz operator spaces.
	If $\psi: \D \to \D$ is a disk automorphism, then one can check that the weighted composition operator
	\begin{equation*}
		U_{\psi}: \mathcal{K}_u \to \mathcal{K}_{u \circ \psi}, \qquad U_{\psi} f = \sqrt{\psi'} (f \circ \psi),
	\end{equation*}
	is unitary.  In particular, this implies that the map
	\begin{equation*}
		\Lambda_{\psi}: \TTO_u \to \TTO_{u \circ \psi}, \qquad \Lambda_{\psi}(A) = U_{\psi} A U_{\psi}^{*},
	\end{equation*}
	which satisfies the useful relationship $\Lambda_{\psi}(A^{u}_{\phi}) = A^{u \circ \psi}_{\phi \circ \psi}$, implements a spatial
	isomorphism between $\TTO_u$ and $\TTO_{u \circ \psi}$.
	
	Another family of spatial isomorphisms arises from the so-called \emph{Crofoot transforms} \cite{Crofoot}
	(see also \cite[Sect.~13]{Sarason}).  For $a \in \D$ and $\psi_{a} = \frac{z - a}{1 - \overline{a} z}$,
	one can verify that the operator
	\begin{equation*}
		U_a: \mathcal{K}_{u}  \to \mathcal{K}_{\psi_a \circ u}, \qquad U_{a} f = \frac{\sqrt{1 - |a|^2}}{1 - \overline{a} u} f,
	\end{equation*}
	is unitary. In particular, the corresponding map
	\begin{equation*}
		\Lambda_a: \TTO_u \to \TTO_{\psi_a \circ u}, \qquad \Lambda_a (A) = U_{a} A U_{a}^{*},
	\end{equation*}
	implements a spatial isomorphism between $\TTO_u$ and $\TTO_{\psi_a \circ u}$.

	Finally, let us define
	\begin{equation*}
		[U_{\#} f](\zeta) = \overline{\zeta} f(\overline{\zeta}) u^{\#}(\zeta),
		\qquad u^{\#}(z) := \overline{u(\overline{z})}.
	\end{equation*}
	The operator $U_{\#}: \mathcal{K}_u \to \mathcal{K}_{u^{\#}}$ is unitary and if
	\begin{equation*}
		\Lambda_{\#}: \TTO_u \to \TTO_{u^{\#}}, \qquad \Lambda_{\#}(A) = U_{\#} A U_{\#}^{*},
	\end{equation*}
	then $\Lambda_{\#}(A^{u}_{\phi}) = A^{u^{\#}}_{\overline{\phi^{\#}}}$
	whence $\TTO_u$ is spatially isomorphic to $\TTO_{u^{\#}}$. Needless to say,
	the three classes of unitary operators $U_{\psi}$, $U_a$, and $U_{\#}$ introduced above should not be confused with
	the Clark operators $U_{\alpha}$ \eqref{eq-ClarkOperator}, which play no role here.

	It turns out that any spatial isomorphism between truncated Toeplitz operator spaces
	can be written in terms of the three basic types described above \cite{TTOSIUES}.

	\begin{theorem}[Cima-Garcia-Ross-Wogen]
		For two inner functions $u_1$ and $u_2$ the spaces $\TTO_{u_1}$ and $\TTO_{u_2}$
		are spatially isomorphic if and only if either $u_1 = \psi \circ u_2 \circ \phi$ or
		$u_1 = \psi \circ u_{2}^{\#} \circ \phi$ for some disk automorphisms $\phi, \psi$.
		Moreover, any spatial isomorphism $\Lambda: \TTO_{u_1} \to \TTO_{u_2}$ can
		be written as $\Lambda = \Lambda_a \Lambda_{\psi}$ or
		$\Lambda_a \Lambda_{\#} \Lambda_{\psi}$, where we allow $a = 0$ or $\psi(z) = z$.
	\end{theorem}
	
	The preceding theorem leads immediately to the following question.
	
	\begin{Question}
		Determine practical conditions on inner functions $u_1$ and $u_2$ which ensure that
		$u_1 = \psi \circ u_2 \circ \phi$ or $u_1 = \psi \circ u_{2}^{\#} \circ \phi$ for some disk automorphisms $\phi, \psi$.
		For instance, do this when $u_1$ and $u_2$ are finite Blaschke products having the same number of zeros,
		counted according to multiplicity.
	\end{Question}

	In the case where one of the inner functions is $z^n$, there is a complete answer \cite{TTOSIUES}.

	\begin{corollary}
		For a finite Blaschke product $u$ of order $n$,
		$\TTO_u$ is spatially isomorphic to $\TTO_{z^n}$ if and only if either
		$u$ has one zero of order $n$ or $u$ has $n$ distinct zeros all lying
		on a circle $\Gamma$ in $\mathbb{D}$ with the property that if these zeros
		are ordered according to increasing argument on $\Gamma$, then
		adjacent zeros are equidistant in the hyperbolic metric.
	\end{corollary}

\section{Algebras of Truncated Toeplitz Operators}\label{SectionAlgebras}
	Recall that $\TTO_u$ is a weakly closed subspace of $\B(\h)$ (see Section \ref{SectionLinearSpace}).
	Although $\TTO_u$ is not an algebra, there are many interesting algebras
	contained within $\TTO_u$.  In fact, the recent thesis \cite{SedlockThesis} and subsequent
	paper of Sedlock \cite{Sed} described them all.  We discuss
	the properties of these so-called \emph{Sedlock algebras} below, along with several
	further results from \cite{GRW}.
		
	To begin with, we require the following generalization (see \cite[Sect.~10]{Sarason})
	of the Clark unitary operators \eqref{eq-ClarkOperator}:
	\begin{equation}\label{eq-Salpha}
		S_{u}^{a} = A_{z}^{u} + \frac{a}{1 - \overline{u(0)} a} k_{0} \otimes C k_0,
	\end{equation}
	where the parameter $a$ is permitted to vary over the closed unit disk $\D^{-}$ (we prefer to
	reserve the symbol $\alpha$ to denote complex numbers of unit modulus).  The operators $S_{u}^{a}$
	turn out to be fundamental to the study of Sedlock algebras.  Before proceeding, let us recall a few basic definitions.
	
	For $A \in \TTO_u$, the commutant $\{A\}'$ of $A$ is defined to be the set of all
	bounded operators on $\K_u$ which commute with $A$.   The weakly closed linear span of
	$\{A^n: n \geq 0\}$ will be denoted by $\mathcal{W}(A)$. Elementary operator theory says that
	$\mathcal{W}(A) \subseteq \{A\}'$ holds and that $\{A\}'$ is a weakly closed subset of $\B(\K_u)$.
	The relevance of these concepts lies in the following two results from \cite[p.~515]{Sarason}
	and \cite{GRW}, respectively.

	\begin{theorem}[Sarason]\label{TheoremSarasonTTOCommute}
		For each $a \in \D^{-}$, $\{S_{u}^{a}\}' \subseteq \TTO_u$.
	\end{theorem}

	\begin{theorem}[Garcia-Ross-Wogen]
		For each $a \in \D^{-}$, $\{S_{u}^{a}\}' = \mathcal{W}(S_{u}^{a})$.
	\end{theorem}
	
	The preceding two theorems tell us that  $\mathcal{W}(S_{u}^{a})$ and
	$\mathcal{W}((S_{u}^{b})^{*})$, where $a, b$ belong to $\D^{-}$,
	are algebras contained in $\TTO_u$.   We adopt the following notation introduced by Sedlock \cite{Sed}:
	\begin{equation*}
		\mathscr{B}_{u}^{a} :=
		\begin{cases} \mathcal{W}(S_{u}^{a}) &\mbox{if } a \in \D^{-}, \\[3pt]
			\mathcal{W}((S_{u}^{1/\overline{a}})^{*}) & \mbox{if } a \in \widehat{\C} \setminus \D^{-}.
		\end{cases}
	\end{equation*}
	Note that $\mathscr{B}_{u}^{0}$ is the algebra of analytic truncated Toeplitz operators
	(i.e., $\mathscr{B}_{u}^{0} = \mathcal{W}(A^{u}_{z})$) and that $\mathscr{B}_{u}^{\infty}$ is the algebra
	of co-analytic truncated Toeplitz operators (i.e., $\mathscr{B}_{u}^{\infty} = \mathcal{W}(A^{u}_{\overline{z}})$).
	The following theorem of Sedlock asserts that the algebras $\mathscr{B}_{u}^{a}$ for $a \in \widehat{\C}$
	are the only maximal abelian algebras in $\TTO_u$.
	
	\begin{theorem}[Sedlock \cite{Sed}]
		If $A, B \in \TTO_u \setminus \{\C I, 0\}$, then $A B \in \TTO_u$ if and only if
		$A, B \in \mathscr{B}_{u}^{a}$ for some $a \in \widehat{\C}$. Consequently,
		every weakly closed algebra in $\TTO_u$ is abelian and is contained in
		some $\mathscr{B}_{u}^{a}$.
	\end{theorem}
	
	Let us gather together a few facts about the Sedlock algebras $\mathscr{B}_{u}^{a}$,
	all of which can be found in Sedlock's paper \cite{Sed}. First we note that
	\begin{equation*}
		(\mathscr{B}_{u}^{a})^{*} = \mathscr{B}_{u}^{1/\overline{a}},
	\end{equation*}
	and
	\begin{equation*}\qquad\qquad
		\mathscr{B}_{u}^{a} \cap \mathscr{B}_{u}^{b} = \C I, \qquad a \not = b.
	\end{equation*}
	Most importantly, we have the following concrete description of $\mathscr{B}_{u}^{a}$.

	\begin{theorem}[Sedlock \cite{Sed}] \label{Hi}
		If $a \in \D$, then
		\begin{equation*}
			\mathscr{B}_{u}^{a} = \left\{A^{u}_{\frac{\phi}{1 - a \overline{u}}}: \phi \in H^{\infty}\right\}.
		\end{equation*}
		Furthermore, if $\phi, \psi \in H^{\infty}$, then we have the following product formula
		\begin{equation*}
			A^{u}_{\frac{\phi}{1 - a \overline{u}}} A^{u}_{\frac{\psi}{1 - a \overline{u}}}
			= A^{u}_{\frac{\phi \psi}{1 - a \overline{u}}}.
		\end{equation*}
	\end{theorem}

	In particular, if $a \in \widehat{\C} \setminus \dD$, then every operator in $\mathscr{B}_{u}^{a}$
	is a truncated Toeplitz operator which can be represented using a bounded symbol.  On the other
	hand, if $a$ belongs to $\dD$, then there may be operators in $\mathscr{B}_{u}^{a}$ which do not have bounded
	symbols. In fact, if $u$ has an ADC at some $\zeta$ in $\dD$, then $k_{\zeta} \otimes k_{\zeta}$
	belongs to $\mathscr{B}_{u}^{u(\zeta)}$.  From here, one can use the example from the remarks after
	Theorem \ref{no-bd-symbol} to produce an operator in $\mathscr{B}_{u}^{u(\zeta)}$ which
	has no bounded symbol.

	Let us now make a few remarks about normal truncated Toeplitz operators.
	For $a$ in $\dD$ the Sedlock algebra $\mathscr{B}_{u}^{a}$ is generated by a unitary operator
	(i.e., a Clark operator) and is therefore an abelian algebra of normal operators. When
	$a \in \widehat{\C} \setminus \dD$, the situation drastically changes \cite{GRW}.

	\begin{theorem}[Garcia-Ross-Wogen]
		If $a \in \widehat{\C} \setminus \dD$, then $A \in \mathscr{B}_{u}^{a}$
		is normal if and only if $A \in \C I$.
	\end{theorem}
	
	In Section \ref{SectionSpatialIsomorphism}, we characterized all possible spatial isomorphisms
	between truncated Toeplitz operator spaces.  In particular, recall that the basic spatial isomorphisms
	$\Lambda_{\psi}, \Lambda_{\#}, \Lambda_{a}$ played a key role.  Let us examine their effect on
	Sedlock algebras.  The following result is from \cite{GRW}.
	
	\begin{theorem}[Garcia-Ross-Wogen]
		For $u$ inner, $a \in \widehat{\C}$, and $c \in \D$, we have
		\begin{equation*}
			\Lambda_{\psi}(\mathscr{B}_{u}^{a}) = \mathscr{B}_{u \circ \psi}^{a}, \qquad
			\Lambda_{\#}(\mathscr{B}_{u}^{a}) = \mathscr{B}_{u^{\#}}^{1/a}, \qquad
			\Lambda_{c}(\mathscr{B}_{u}^{a}) = \mathscr{B}_{u_c}^{\ell_{c}(a)},
		\end{equation*}
		where
		\begin{equation*}
			u_{c} = \frac{u - c}{1 - \overline{c} u}, \qquad \ell_{c}(a) =
			\begin{cases}
				\frac{a - c}{1 - \overline{c} a} &\mbox{if } a \not = \frac{1}{\overline{c}}, \\[3pt]
				\infty & \mbox{if } a = \frac{1}{\overline{c}}.
			\end{cases}
		\end{equation*}
	\end{theorem}
	
	For $a$ in $\widehat{\C} \setminus \dD$, the preceding theorem follows from direct computations based
	upon Theorem \ref{Hi}.  When $a$ belongs to $\dD$, however, a different proof is required. One can
	go even further and investigate when two Sedlock algebras are
	spatially isomorphic.  The following results are from  \cite{GRW}.
	
	\begin{theorem}[Garcia-Ross-Wogen]
		If $u(z) = z^n$ and $a, a' \in \D$, then
		$\mathscr{B}_{u}^{a}$ is spatially isomorphic to $\mathscr{B}_{u}^{a'}$ if and only if $|a| = |a'|$.
	\end{theorem}
	
	\begin{theorem}[Garcia-Ross-Wogen]
		If
		\begin{equation*}
			u(z)= \exp\left(\frac{z + 1}{z - 1}\right)
		\end{equation*}
		and $a, a' \in \D$, then $\mathscr{B}_{u}^{a}$ is spatially isomorphic to $\mathscr{B}_{u}^{a'}$ if and only if $a = a'$.
	\end{theorem}

	Before moving on, let us take a moment to highlight an interesting operator integral formula
	from \cite[Sect.~12]{Sarason} which is of some relevance here.
	Recall that for $\alpha$ in $\dD$, the operator $S_{u}^{\alpha}$ given by \ref{eq-Salpha} is unitary, whence
	$\phi(S_{u}^{\alpha})$ is defined by the functional calculus for $\phi$ in $L^{\infty}$.  By Theorem
	\ref{TheoremSarasonTTOCommute}
	we see that $\phi(S_{u}^{\alpha})$ belongs to $\TTO_u$.   Using the Aleksandrov
	disintegration theorem (Theorem \ref{Alek-Dec-Form}), one can then prove that
	
		\begin{equation*}\qquad\qquad\qquad
		\langle A_{\phi} f, g \rangle = \int_{\dD} \langle \phi(S_{u}^{\alpha}) f, g \rangle
		\,dm(\alpha), \qquad f, g \in \mathcal{K}_{u},
	\end{equation*}
	which can be written in the more compact and pleasing form
	\begin{equation*}
		A_{\phi} = \int_{\dD} \phi(S_{u}^{\alpha}) \,dm(\alpha).
	\end{equation*}
	A similar formula exists for symbols $\phi$ in $L^2$, but the preceding formulae
	must be interpreted carefully since the operators $\phi(S_{u}^{\alpha})$ may be unbounded.

	Recall that for each $\phi$ in $L^{\infty}$, the Ces\`aro means of $\phi$ are
	trigonometric polynomials $\phi_n$ which approximate $\phi$ in the weak-$*$ topology
	of $L^{\infty}$.   From the discussion above and Corollary \ref{wkly-dense} we know that
	\begin{equation*}
		\{q(S_{u}^{\alpha}): \mbox{$q$ is a trigonometric polynomial}, \alpha \in \dD\}
	\end{equation*}
	is weakly dense in $\TTO_u$. When $u$ is a finite Blaschke product, it turns out that we can do much better.
	The following result can be found in \cite{MR2440673}, although it can be gleaned from
	\cite[Sect.~12]{Sarason}.

	\begin{theorem}
		Let $u$ be a Blaschke product of degree $N$ and let $\alpha_1, \alpha_2 \in \dD$ with
		$\alpha_1 \not = \alpha_2$. Then for any $\varphi \in L^2$, there are polynomials $p, q$
		of degree at most $N$ so that $A_{\varphi} = p(S_{u}^{\alpha_1}) + q(S_{u}^{\alpha_2})$.
	\end{theorem}

\section{Truncated Toeplitz $C^{*}$-Algebras}\label{SectionCStar}

	In the following, we let $\h$ denote a separable complex Hilbert space.  For each $\mathscr{X} \subseteq \mathcal{B}(\h)$,
	let $C^*(\mathscr{X})$ denote the unital $C^*$-algebra generated by $\mathscr{X}$.  In other words, $C^*(\mathscr{X})$
	is the closure, in the norm of $\B(\h)$, of the unital algebra generated by the operators in $\mathscr{X}$ and their adjoints.
	Since we are frequently
	interested in the case where $\mathscr{X} = \{A\}$ is a singleton, we often write $C^*(A)$ in place of $C^*(\{A\})$
	in order to simplify our notation.
	
	Recall that the \emph{commutator ideal} $\mathscr{C}(C^*(\mathscr{X}))$
	of $C^*(\mathscr{X})$, is the smallest closed two-sided ideal which contains the \emph{commutators}
	\begin{equation*}
		[A,B]:=AB-BA
	\end{equation*}
	where $A$ and $B$ range over all elements of $C^*(\mathscr{X})$.  Since
	the quotient algebra $C^*(\mathscr{X})/ \mathscr{C}(C^*(\mathscr{X}))$ is an abelian $C^*$-algebra,
	it is isometrically $*$-isomorphic to $C(Y)$, the set of all continuous functions on some compact Hausdorff
	space $Y$ \cite[Thm.~1.2.1]{MR1721402}.  We denote this relationship
	\begin{equation}\label{eq-CCCY}
		\frac{C^*(\mathscr{X})}{\mathscr{C}(C^*(\mathscr{X}))} \cong C(Y).
	\end{equation}
	Putting this all together, we have the short exact sequence
	\begin{equation}\label{eq-ShortExact}
		0 \longrightarrow \mathscr{C}(C^*(\mathscr{X})) \overset{\iota}{\longrightarrow}
		C^*(\mathscr{X}) \overset{\pi}{\longrightarrow} C(Y) \longrightarrow 0,
	\end{equation}
	where $\iota: \mathscr{C}(C^*(\mathscr{X})) \to C^*(\mathscr{X})$ is the inclusion map
	and $\pi:C^*(\mathscr{X})\to C(Y)$ is the composition of the quotient map with
	the isometric $*$-isomorphism which implements \eqref{eq-CCCY}.

	The \emph{Toeplitz algebra} $C^*(T_z)$, where $T_z$ denotes the unilateral shift on the classical Hardy space $H^2$,
	has been extensively studied since the seminal work of Coburn in the late 1960s \cite{Coburn, Coburn2}.
	Indeed, the Toeplitz algebra is now
	one of the standard examples discussed in many well-known texts (e.g., \cite[Sect.~4.3]{MR1865513},
	\cite[Ch.~V.1]{Davidson}, \cite[Ch.~7]{MR1634900}).  In this setting, we have $\mathscr{C}(C^*(T_z)) = \mathscr{K}$,
	the ideal of compact operators on $H^2$, and $Y = \dD$, so that the short exact sequence
	\eqref{eq-ShortExact} takes the form
	\begin{equation}\label{eq-CoburnExact}
		0 \longrightarrow \mathscr{K} \overset{\iota}{\longrightarrow} C^*(T_z)
		\overset{\pi}{\longrightarrow} C(\dD) \longrightarrow 0.
	\end{equation}
	In other words, $C^*(T_z)$ is an \emph{extension} of $\mathscr{K}$ by $C(\dD)$.
	It also follows that
	\begin{equation*}
		C^*(T_z) = \{ T_{\phi} + K : \phi \in C(\dD), K \in \mathscr{K}\},
	\end{equation*}
	and, moreover, that each element of $C^*(T_z)$ enjoys a unique decomposition of the form
	$T_{\phi} + K$ \cite[Thm.~4.3.2]{MR1865513}.  As a consequence, we see that the map $\pi:C^*(T_z)\to C(\dD)$ is given
	by $\pi(T_{\phi}+K) = \phi$.
	
	Needless to say, the preceding results have spawned numerous generalizations and variants over the years.
	For instance, one can consider $C^*$-algebras generated by matrix-valued Toeplitz operators
	or Toeplitz operators acting on other Hilbert function spaces (e.g., the Bergman space \cite{ACM}).
	As another example, if $\mathscr{X}$ denotes the space of functions on $\dD$ which are
	both piecewise and left continuous, then Gohberg and Krupnik
	proved that $\mathscr{C}(C^*(\mathscr{X})) = \mathscr{K}$ and obtained the short exact sequence
	\begin{equation*}
		0 \longrightarrow \mathscr{K} \overset{\iota}{\longrightarrow}
		C^*(\mathscr{X}) \overset{\pi}{\longrightarrow} C(Y) \longrightarrow 0,
	\end{equation*}
	where $Y$ is the cylinder $\dD \times [0, 1]$, endowed with a nonstandard topology \cite{GK}.

	In the direction of truncated Toeplitz operators, we have the following analogue of Coburn's work.
		
	\begin{theorem} \label{TheoremMainContinuous}
		If $u$ is an inner function, then
		\begin{enumerate}\addtolength{\itemsep}{0.5\baselineskip}
			\item $\mathscr{C}(C^*(A_z^u)) = \mathscr{K}^u$, the algebra of compact operators on $\mathcal{K}_{u}$,
			\item $C^*(A_z^u)/\mathscr{K}^u$ is isometrically $*$-isomorphic to $C(\sigma(u) \cap \dD)$,
			\item If $\phi \in C(\dD)$, then $A_{\phi}^u$ is compact if and only if $\phi(\sigma(u)\cap \dD) = \{0\}$,
			\item $C^*(A_z^u) = \{ A_{\phi}^u  + K : \phi \in C(\dD), K \in \mathscr{K}^u\}$,
			%\item Each element of $C^*(A_z^u)$ has a representation of the form $A_{\phi}^u + K$, where $\phi \in C(\dD)$ and $K \in \mathscr{K}^{u}$,
			\item If $\phi \in C(\dD)$, then $\sigma_{e}(A^{u}_{\phi}) = \phi(\sigma_{e}(A_z^u))$,
			\item For $\phi \in C(\dD)$, $\|A^{u}_{\phi}\|_{e}  = \sup\{ |\phi(\zeta)| : \zeta \in \sigma(u)\cap\dD\}$.
		\end{enumerate}
		%In particular,
		%\begin{equation*}
		%	0\longrightarrow\mathscr{K}^u \overset{\iota}{\longrightarrow} C^*(A_z^u)
		%	\overset{\pi}{\longrightarrow} C(\sigma(u) \cap \dD) \longrightarrow 0
		%\end{equation*}
		%is a short exact sequence, where $\pi$ is the composition of the quotient map $C^{*}(A^{u}_{z}) \to C^{*}(A^{u}_{z})/\mathscr{K}^{u}$ together with the map implementing the $*$-isomorphism $C^{*}(A^{u}_{z})/\mathscr{K}^{u}  \cong C(\sigma(u) \cap \dD)$. Thus $C^*(A_z^u)$ is an extension of the compact operators
		%by $C(\sigma(u) \cap \dD)$.
			\end{theorem}

	In recent work \cite{GWR-Cstar}, the authors and W.~Wogen were able to provide operator
	algebraic proofs of the preceding results, utilizing an approach similar in spirit to the original work of Coburn.
	However, it should also be noted that many of the statements in Theorem \ref{TheoremMainContinuous}
	can be obtained using the explicit triangularization theory developed by Ahern and Clark in \cite{AC70a}
	(see the exposition in \cite[Lec.~V]{N1}).

\section{Unitary Equivalence to a Truncated Toeplitz Operator}\label{SectionCSO}
	A significant amount of evidence is mounting that truncated Toeplitz operators may play a significant
	role in some sort of model theory for complex symmetric operators \cite{TTOSIUES, ATTO,STZ}.  At this point,
	however, it is still too early to tell what exact form such a model theory should take.
	On the other hand, a surprising array of complex symmetric
	operators can be concretely realized in terms of truncated Toeplitz operators (or direct sums of such operators),
	without yet even venturing to discuss vector-valued truncated Toeplitz operators.

	Before discussing unitary equivalence, however, we should perhaps say a few words about
	similarity.  A number of years ago, D.S.~Mackey, N.~Mackey, and Petrovic asked whether or not the
	inverse Jordan structure problem can be solved
	in the class of Toeplitz matrices \cite{Mackey}.  In other words, given any Jordan canonical form, can one find a
	Toeplitz matrix which is similar to this form? A negative answer to this question was subsequently provided by Heinig
	in \cite{Heinig}. On the other hand, it turns out that the inverse Jordan structure problem is always solvable
	in the class of truncated Toeplitz operators, for we have the following theorem \cite[Thm.~6.2]{TTOSIUES}.
	
	\begin{theorem}[Cima-Garcia-Ross-Wogen]
		Every operator on a finite dimensional space is similar to a co-analytic truncated Toeplitz operator.
	\end{theorem}

	In light of the preceding theorem, it is clear that simple, purely algebraic,  tools will be
	insufficient to settle the question of whether every complex symmetric operator can be represented
	in some fashion using truncated Toeplitz operators.  We turn our attention now toward unitary equivalence.
	
	Let us begin by recalling an early result of Sarason, who observed that the Volterra integration operator
	\eqref{eq-Volterra},
	a standard example of a complex symmetric operator \cite{G-P, G-P-II, CCO},
	is unitarily equivalent to a truncated Toeplitz operator acting on the $\K_u$
	space corresponding to the atomic inner function $u(z) = \exp( \frac{z+1}{z-1})$ \cite{MR0192355} (although the term
	``truncated Toeplitz operator'' was not yet coined).
	Detailed computations using the theory of model operators and characteristic functions
	can be found in \cite[p.~41]{N3}.
	
	What was at first only an isolated result has recently begun to be viewed as a seminal observation.  More recently,
	a number of standard classes of complex symmetric operators have been identified as being unitarily
	equivalent to truncated Toeplitz operators.  Among the first observed examples are
	\begin{enumerate}\addtolength{\itemsep}{0.5\baselineskip}
		\item rank-one operators \cite[Thm.~5.1]{TTOSIUES},
		\item $2 \times 2$ matrices \cite[Thm.~5.2]{TTOSIUES},
		\item normal operators \cite[Thm.~5.6]{TTOSIUES},
		\item for $k \in \N \cup \{\infty\}$, the $k$-fold inflation of a finite Toeplitz matrix \cite[Thm.~5.7]{TTOSIUES}
	\end{enumerate}
	This last item was greatly generalized by Strouse, Timotin, and Zarrabi \cite{STZ}, who proved that a remarkable
	array of inflations of truncated Toeplitz operators are themselves truncated Toeplitz operators.  In addition,
	a variety of related results concerning tensor products, inflations, and direct sums are given in \cite{STZ}.
	The key to many of these results lies in the fact that if $B$ is an inner function, then
	\begin{equation*}
		h \otimes f \mapsto h(f\circ B)
	\end{equation*}
	extends linearly to a unitary operator $\Omega_B: \K_B\otimes L^2 \to L^2$ and, moreover,
	this operator maps $\K_B \otimes H^2$ onto $H^2$.  Letting $\omega_B: \K_B\otimes \K_u \to \K_{u\otimes B}$
	denote the restriction of $\Omega_B$ to $\K_B\otimes \K_u$, one can obtain the following general theorem.

	\begin{theorem}[Strouse-Timotin-Zarrabi \cite{STZ}]
		Let $B$ and $u$ be inner functions, and suppose that $\psi,\phi$ belong to $L^2$
		and satisfy the following conditions:
		\begin{enumerate}\addtolength{\itemsep}{0.5\baselineskip}
			\item The operators $A_{\overline{B}^j \psi}^B$ are bounded,
				and nonzero only for a finite number of $j \in \Z$.
			\item $A_{\phi}^u$ is bounded.
			\item $\psi(\phi\circ B) \in L^2$.
		\end{enumerate}
		Then $A_{\psi(\phi\circ B)}^{u\circ B}$ is bounded and
		\begin{equation*}
			A_{\psi(\phi\circ B)}^{u\circ B} \omega_B = \omega_B \left( \sum_j (A_{\overline{B}^j \psi}^B \otimes A_{z^j \phi}^u) \right).
		\end{equation*}
	\end{theorem}
		
	We state explicitly only a few more results from \cite{STZ}, hoping to give the reader the general flavor of this surprising work.
	In the following, we say that the inner function $B$ is of order $n$ if $B$ is a finite Blaschke product of degree $n$,
	and of order infinity otherwise.

	\begin{theorem}[Strouse-Timotin-Zarrabi]
		Suppose that $u$ is an inner function, that $\phi \in L^2$, and that $B$ is
		an inner function of order $k$ for some $k \in \N\cup\{\infty\}$.  Assume also that $A_{\phi}^u$ is
		bounded. Then $A_{\phi\circ B}^{u\circ B}$ is bounded and unitarily equivalent to $I_k \otimes A_{\phi}^{u}$.
	\end{theorem}
	
	\begin{theorem}[Strouse-Timotin-Zarrabi]
		If $\psi$ is an analytic function, $A_{\psi}^B$ is bounded, and $R$ is a non-selfadjoint operator of rank one,
		then $A_{\psi}^B \otimes R$ is unitarily equivalent to a truncated Toeplitz operator.
	\end{theorem}
	
	\begin{theorem}[Strouse-Timotin-Zarrabi]
		Suppose $u$ is inner, $\phi \in H^{\infty}$, and $(A_{\phi}^u)^2 = 0$.
		If $k = \dim \K_u \ominus \ker A_{\phi}^u$, then $A_{\phi}^u \oplus 0_k$ is unitarily equivalent to a truncated
		Toeplitz operator.
	\end{theorem}

	Although a few results concerning matrix representations of
	truncated Toeplitz operators have been obtained \cite{MR2440673, TTOSIUES, STZ}, the general question of
	determining whether a given matrix represents a truncated Toeplitz operator, with respect to some orthonormal basis
	of some $\K_u$ space,
	appears difficult.   On the other hand, it is known that every truncated Toeplitz operator
	is unitarily equivalent to a complex symmetric matrix \cite{G-P,CCO}, a somewhat more general issue
	which has been studied for its own independent interest \cite{UECSMGC, UECSMMC, UECSMLD, Tener, Vermeer}.
	
	The main result of \cite{ATTO} is the following simple criterion  for determining
	whether or not a given matrix is unitarily equivalent to a trunctated
	Toeplitz operator having an analytic symbol.

	\begin{theorem}[Garcia-Poore-Ross]\label{TheoremAnalytic}
		Suppose  $M \in {\bf M}_n(\C)$ has distinct eigenvalues $\lambda_1, \lambda_2, \ldots, \lambda_n$ with corresponding
		unit eigenvectors $\vec{x}_1,\vec{x}_2,\ldots,\vec{x}_n$. Then
		$M$ is unitarily equivalent to an analytic truncated Toeplitz operator, on some model space $\mathcal{K}_u$
		if and only if there are distinct points $z_1,z_2, \ldots, z_{n-1}$ in $\D$ such that
		\begin{equation}\label{eq-TripleExplicit}
			\inner{\vec{x}_n,\vec{x}_i}\inner{\vec{x}_i,\vec{x}_j} \inner{\vec{x}_j,\vec{x}_n}
			= \frac{ (1 - |z_i|^2)(1-|z_j|^2) }{ 1 - \overline{z_j}z_i}
		\end{equation}
		holds for $1 \leqslant i \leqslant j <n$.
	\end{theorem}
	
	The method of Theorem \ref{TheoremAnalytic} is constructive, in the sense that if \eqref{eq-TripleExplicit}
	is satisfied, then one can construct an inner function $u$ and a polynomial $\phi$ such that
	$M$ is unitarily equivalent to $A_{\phi}^{u}$.
	In fact, $u$ is the Blaschke product having zeros at $z_1,z_2,\ldots,z_{n-1}$ and $z_n = 0$.
	Using Theorem \ref{TheoremAnalytic} and other tools, one can prove the following result from \cite{ATTO}.
	
	\begin{theorem}[Garcia-Poore-Ross]
		Every complex symmetric operator on a $3$-dimensional Hilbert space is unitarily
		equivalent to a direct sum of truncated Toeplitz operators.
	\end{theorem}

	Taken together, these results from \cite{TTOSIUES, ATTO, STZ} yield a host of open questions,
	many of which are still open, even in the finite dimensional setting. 	
	
	\begin{Question}\label{QuestionMain}
		Is every complex symmetric matrix $M \in {\bf M}_n(\C)$ unitarily equivalent to a
		direct sum of truncated Toeplitz operators?
	\end{Question}	
	
	\begin{Question}
		Let $n\geq 4$.  Is every irreducible complex symmetric matrix $M \in {\bf M}_n(\C)$
		unitarily equivalent to a truncated Toeplitz operator?
	\end{Question}

	Recently, the first author and J.~Tener \cite{UET} showed that every finite complex symmetric
	matrix is unitarily equivalent to a direct sum of (i) irreducible complex symmetric matrices
	or (ii) matrices of the form $A \oplus A^T$ where $A$ is irreducible
	and not unitarily equivalent to a complex symmetric matrix
	(such matrices are necessarily $6 \times 6$ or larger).  This immediately
	suggests the following question.

	\begin{Question}\label{QuestionTranspose}
		For $A \in {\bf M}_n(\C)$, is
		the matrix $A \oplus A^T \in {\bf M}_{2n}(\C)$ unitarily equivalent to a direct sum of truncated
		Toeplitz operators?
	\end{Question}
	
	One method for producing complex symmetric matrix representations of a given truncated Toeplitz
	operator is to use \emph{modified Clark bases} \eqref{basis-modClark} for $\K_u$.

	\begin{Question}
		Suppose that $M \in {\bf M}_n(\C)$ is complex symmetric.  If $M$
		is unitarily equivalent to a truncated Toeplitz operator, does there exist an inner function $u$,
		a symbol $\phi \in L^{\infty}$, and a modified Clark basis for $\K_u$
		such that $M$ is the matrix representation of $A_{\phi}^{\Theta}$ with respect to this basis?
		In other words, do all such unitary equivalences between complex symmetric matrices and truncated
		Toeplitz operators arise from Clark representations?
	\end{Question}

\section{Unbounded Truncated Toeplitz Operators}\label{SectionUnbounded}
	As we mentioned earlier (Subsection \ref{SubsectionTTO}), for a symbol $\phi$ in $L^2$ and an
	inner function $u$, the truncated Toeplitz operator $A^u_{\phi}$ is closed and densely defined on
	the domain
	\begin{equation*}
		\mathscr{D}(A^u_{\phi}) = \{f \in \mathcal{K}_u: P_u(\phi f) \in \K_{u}\}
	\end{equation*}
	in $\K_u$.  In particular, the analytic function $P_u (\phi f)$ can be defined on $\D$
	by writing the formula \eqref{eq-Pu-formula} as an integral.  In general, we actually have
	$CA^u_{\phi} C = A^u_{\overline{\phi}}$, and, when $A^u_{\phi}$ is bounded (i.e., $A^{u}_{\phi} \in \TTO_u$),
	we have $A^u_{\overline{\phi}} = (A^u_{\phi})^{*}$. Let us also recall an old and important
	result of Sarason \cite{Sarason-NF} which inspired the comuntant lifting theorem \cite{MR2760647}.

	\begin{theorem}[Sarason] \label{Sarason-Lifting-Thm}
		The  bounded operators on $\mathcal{K}_{u}$ which commute with
		$A^{u}_{\overline{z}}$ are $\{A^u_{\overline{\phi}}: \phi \in H^{\infty}\}$.
	\end{theorem}

	In the recent papers \cite{MR2418122, MR2468883}, Sarason studied unbounded
	Toeplitz operators (recall that a Toeplitz operator on $H^2$ is bounded if and only if the symbol is bounded)
	as well as unbounded truncated Toeplitz operators. We give a brief survey of these results.

	In the following, $N$ denotes the \emph{Nevanlinna class}, the set of all quotients $f/g$ where $f$
	and $g$ belong to $H^{\infty}$ and $g$ is non-vanishing on $\D$. The \emph{Smirnov class} $N^{+}$ denotes the
	subset of $N$ for which the denominator $g$ is not only non-vanishing on $\D$ but outer.  By \cite{MR2418122} each $\phi$ in $N^{+}$ can be written uniquely as
	$\phi = \frac{b}{a}$ where $a$ and $b$ belong to $H^{\infty}$, $a$ is outer, $a(0) > 0$, and $|a|^2 + |b|^2 = 1$
	almost everywhere on $\dD$. Sarason calls this the {\em canonical representation} of $\phi$.

	For $\phi$ in $N^{+}$ define the Toeplitz operator $T_{\phi}$ as multiplication by $\phi$
	on its domain $\mathscr{D}(T_{\phi}) = \{f \in H^2: \phi f \in H^2\}$. In particular, observe that there
	is no projection involved in the preceding definition.

	\begin{theorem}[Sarason \cite{MR2418122}]
		For $\phi = b/a \in N^{+}$, written in canonical form, $T_{\phi}$ is a closed operator on
		$H^2$ with dense domain $\mathscr{D}(T_{\phi}) = a H^2$.
	\end{theorem}

	There is no obvious way to define the co-analytic Toeplitz operator $T_{\overline{\phi}}$
	on $H^2$ for $\phi \in N^{+}$. Of course we can always \emph{define} $T_{\overline{\phi}}$
	to be $T_{\phi}^{*}$ and this makes sense when $\phi$ belongs to $H^{\infty}$. In order to
	justify the definition $T_{\overline{\phi}} := T_{\phi}^{*}$ for $\phi \in N^{+}$, however,
	we need to take care of some technical details.

	As the preceding theorem shows, if $\phi \in N^{+}$, then $T_{\phi}$ is a closed operator
	with dense domain $a H^2$.  Basic functional analysis tells us that its adjoint $T_{\phi}^{*}$
	is also closed and densely defined. In fact, one can show that $\mathscr{D}(T_{\phi}^{*} )$
	is the associated deBranges-Rovnyak space $\mathcal{H}(b)$
	 \cite{MR2418122}. In order to understand $T_{\overline{\phi}}$ we proceed,
	at least formally, as we do when examining $T_{\overline{\phi}}$ when $\phi$ is bounded.
	Let $\phi$ and $f$ have Fourier expansions
	\begin{equation*}
		\phi \sim \sum_{n = 0}^{\infty}  \phi_n \zeta^n, \qquad f \sim \sum_{n = 0}^{\infty} f_n \zeta^n.
	\end{equation*}
	Formal series manipulations show that
	\begin{align*}
		T_{\overline{\phi}} f & = P_{+}(\overline{\phi} f)\\
		& = P_{+}\left(\left(\sum_{n = 0}^{\infty} \overline{\phi_n} \overline{\zeta}^n\right)
		 \left(\sum_{m = 0}^{\infty} f_m \zeta^m\right)\right)\\
		& = P_{+}\left(\sum_{n, m = 0}^{\infty} \zeta^{m - n} \overline{\phi_n} f_m \right)\\
		& = P_{+}\left(\sum_{k = -\infty}^{\infty} \zeta^k \sum_{m = 0}^{\infty} \overline{\phi_m} f_{k + m} \right)\\
		& = \sum_{k = 0}^{\infty} \zeta^k \sum_{m = 0}^{\infty} \overline{\phi_m} f_{k + m}.
	\end{align*}
	This suggests that if $\phi = \sum_{n = 0}^{\infty} \phi_n z^n$ is the power series representation
	for $\phi$ in $N^{+}$, then we should define, for each function $f(z) = \sum_{n = 0}^{\infty} f_n z^n$
	analytic in a neighborhood of $\D^{-}$,
	\begin{equation*}
		(t_{\overline{\phi}} f)(z) := \sum_{k = 0}^{\infty} z^k \left( \sum_{m = 0}^{\infty} \overline{\phi_m} f_{k + m}\right).
	\end{equation*}
	It turns out that $t_{\overline{\phi}} f$, so defined, is an analytic function on $\D$.  The following result
	indicates that this is indeed the correct approach to defining $T_{\phi}^*$.
	
	\begin{theorem}[Sarason  \cite{MR2418122}]
		If $\phi \in N^{+}$, then $t_{\overline{\phi}}$ is closable and $T_{\phi}^{*}$ is its closure.
	\end{theorem}

	In light of the preceding theorem, for $\phi$ in $N^{+}$ we may define $T_{\overline{\phi}}$
	to be $T_{\phi}^{*}$.  More generally, we can define, for each $\phi$ in $N^{+}$ and $u$ inner,
	the truncated Toeplitz operator $A^u_{\overline{\phi}}$ by
	\begin{equation*}
		A^u_{\overline{\phi}} := T_{\overline{\phi}}|_{\mathscr{D}(T_{\overline{\phi}}) \cap \mathcal{K}_{u}}.
	\end{equation*}
	We gather up some results about $A^{u}_{\overline{\phi}}$ from \cite{MR2418122}.

	\begin{theorem}[Sarason]
		If $\phi = b/a \in N^{+}$ is in canonical form and $u$ is inner, then
		\begin{enumerate}\addtolength{\itemsep}{0.5\baselineskip}
			\item $A^{u}_{\overline{\phi}}$ is closed and densely defined.
			\item $A^{u}_{\overline{\phi}}$ is bounded if and only if $\mbox{dist}(b, u H^{\infty}) < 1$.
		\end{enumerate}
	\end{theorem}

	From Theorem \ref{TTOareCSO} we know that $A = CA^*C$ whenever $A \in \TTO_u$.
	Therefore it makes sense for us to define
	\begin{equation*}
		A^{u}_{\phi} := C A_{\overline{\phi}} C
	\end{equation*}
	for $\phi \in N^{+}$.
	It turns out that $\mathscr{D}(A^{u}_{\phi}) = C \mathscr{D}(A^{u}_{\phi})$ and that
	the operator $A^{u}_{\phi}$ is closed and densely defined. Fortunately this definition
	makes sense in terms of adjoints.

	\begin{theorem}[Sarason]
		For inner $u$ and $\phi \in N^{+}$, the operators $A_{\phi}$ and $A_{\overline{\phi}}$ are adjoints of each other.
	\end{theorem}

	What is the analog of Theorem \ref{Sarason-Lifting-Thm} for
	of unbounded truncated Toeplitz operators? In \cite{MR2468883} Sarason showed that
	\begin{equation*}
		A^u_{\overline{\phi}} A^u_{\overline{z}} f = A^u_{\overline{z}} A^{u}_{\overline{\phi}} f,
	\end{equation*}
	holds for $f$ in $\mathscr{D}(A^{u}_{\overline{\phi}})$
	and thus one might be tempted to think that the closed densely defined operators which commute with
	$A^{u}_{\overline{z}}$ are simply $\{A_{\overline{\phi}}: \phi \in N^{+}\}$.
	Unfortunately the situation is more complicated and one needs to define $A^{u}_{\overline{\phi}}$
	for a slightly larger class of symbols than $N^{+}$. Sarason works out the details in \cite{MR2468883}
	and identifies the closed densely defined operators on $\mathcal{K}_u$ which commute with
	$A^{u}_{\overline{z}}$ as the operators $A^{u}_{\overline{\phi}}$ where the symbols $\phi$ come
	from a so-called \emph{local Smirnov class} $N^{+}_u$. The details are somewhat
	technical and so we therefore leave it to the reader to explore this topic
	further in Sarason paper \cite{MR2468883}.

\section{Smoothing Properties of Truncated Toeplitz Operators}
	Let us return to Theorem \ref{AC-paper},  an important
	result of Ahern and Clark which characterizes those functions in the model space
	$\K_u$ which have a finite angular derivative in the sense of Carath\'eodory (ADC) at some point $\zeta$ on $\dD$.
	In particular, recall that every function in $\K_u$ has a finite nontangential limit at $\zeta$ precisely when
	$u$ has an ADC at $\zeta$.  The proof of this ultimately relies on the fact that this statement is equivalent
	to the condition that $(I - \overline{\lambda} A^{u}_{z})^{-1} P_u 1$ is bounded as $\lambda$ approaches
	$\zeta$ nontangentially.  One can see this by observing that
	\begin{equation*}
		\langle f, (I - \overline{\lambda} A^{u}_{z})^{-1} P_u 1\rangle = f(\lambda)
	\end{equation*}
	holds for all $f$ in $\K_u$.  If one replaces $P_u 1$ in the formula above with $P_u h$ for some
	$h$ in $H^{\infty}$, then a routine calculation shows that
	\begin{equation*}
		\langle f, (I - \overline{\lambda} A^{u}_{z})^{-1} P_u h\rangle = (A^{u}_{\overline{h}} f)(\lambda)
	\end{equation*}
	for all $f$ in $\K_u$.  An argument similar to that employed by Ahern and Clark shows that
	$(A^{u}_{\overline{h}} f)(\lambda)$ has a finite nontangential limit at $\zeta$ for each $f$
	in $\mathcal{K}_u$ if and only if $(I - \overline{\lambda} A^{u}_{z})^{-1} P_u h$ is bounded as
	$\lambda$ approaches $\zeta$ nontangentially.
	
	Let us examine the situation when $u$ is an infinite Blaschke product with
	zeros $\{\lambda_{n}\}_{n \geq 1}$, repeated according to multiplicity.
	Recall that the Takenaka basis $\{\gamma_{n}\}_{n \geq 1}$, defined by \eqref{basis-Tak},
	is an orthonormal basis for $\K_u$.  For each $\zeta$ in $\dD$, a calculation
	from \cite{RH} shows that $A^{u}_{\overline{h}} \gamma_n$ is a rational function (and so can be defined at any $\zeta \in \dD$).
	From this one can obtain the following analogue of the Ahern-Clark result \cite{RH}.

	\begin{theorem}[Hartmann-Ross]\label{TheoremSmoothing}
		If $u$ is a Blaschke product with zeros $\{\lambda_{n}\}_{n \geq 1}$ and $h \in H^{\infty}$,
		then every function in $\ran A^{u}_{\overline{h}}$ has a finite nontangential limit at $\zeta \in \dD$ if and only if
		\begin{equation*}
			\sum_{n  = 1}^{\infty} |(A^{u}_{\overline{h}} \gamma_{n})(\zeta)|^2 < \infty.
		\end{equation*}
	\end{theorem}

	From here one can see the smoothing properties of the co-analytic truncated Toeplitz operator $A^{u}_{\overline{h}}$.
	If $u$ happens to be an interpolating Blaschke product, then the condition in the above theorem reduces to
	\begin{equation*}
		\sum_{n = 1}^{\infty} (1 - |\lambda_n|^2) \left| \frac{h(\lambda_n)}{\zeta - \lambda_{n}}\right|^2 < \infty.
	\end{equation*}
	
	The following open problem now suggests itself.
	
	\begin{Question}
		Obtain extensions of Theorem \ref{TheoremSmoothing} to general inner functions $u$
		and symbols $h \in L^2$.
	\end{Question}

\section{Nearly Invariant Subspaces}
	We conclude this survey with a few remarks about truncated Toeplitz operators which act on
	a family of spaces that are closely related to the model spaces $\K_u$.  To be more precise,
	we say that a (norm closed) subspace $\M$ of $H^2$ is \emph{nearly invariant} if the following
	divisibility condition holds
	\begin{equation} \label{nearly-definition}
		f \in \M,\,\, f(0) = 0\quad \Longrightarrow \quad \frac{f}{z} \in \M.
	\end{equation}
	These spaces were first considered and characterized in \cite{Hi88, Sa88} and
	they continue to be the focus of intense study \cite{AR96, HSS04, MP05, KN06, AK08, CCP}.
	
	The link between nearly invariant subspaces of $H^2$ and model spaces
	is supplied by a crucial result of Hitt \cite{Hi88}, which
	asserts that there is a unique solution $g$ to the extremal problem
	\begin{equation} \label{extg}
		 \sup\{\Re g(0):g\in \M, \,\|g\| =1\},
	\end{equation}
	and moreover, that there is an inner function $u$ so that
	\begin{equation*}
		\M = g \mathcal{K}_u
	\end{equation*}
	and such that the map $W_{g}: \mathcal{K}_u \to \M$ defined by
	\begin{equation}\label{eq-IsometricMultiplier}
		W_{g} f = g f
	\end{equation}
	is unitary.  The function $g$ is called the \emph{extremal function} for the
	nearly invariant subspace $\M$.  It is important to observe that since $g$ belongs to $\M = g\K_u$,
	the inner function $u$ must satisfy $u(0) = 0$.  We remark that the content of these observations is nontrivial, for
	the set $g\K_u$, for arbitrary $g$ in $H^2$ and $u$ inner, is not necessarily even a subspace of $H^2$
	since it may fail to be closed.
	
	In the other direction, Sarason showed that if $u$ is an inner function which satisfies
	$u(0) = 0$, then every isometric multiplier from $\mathcal{K}_u$ into $H^2$ takes the form
	\begin{equation} \label{Sar-g}
		 g=\frac{a}{1-ub},
	\end{equation}
	where $a$ and $b$ are in the unit ball of $H^{\infty}$ satisfy $|a|^2+|b|^2=1$ a.e.~on $\dD$  \cite{Sa88}.
	Consequently, one sees that $\M=g \mathcal{K}_u$ is a (closed) nearly invariant subspace of $H^2$ with
	extremal function $g$ as in \eqref{extg}.
		
	The next natural step towards defining truncated Toeplitz operators on the nearly invariant subspace $\M = g \mathcal{K}_u$
	is to understand $P_{\M}$, the orthogonal projection of $L^2$ onto $\M$.  The following lemma from \cite{RH} provides
	an explicit formula relating $P_{\M}$ and $P_u$.

	\begin{lemma}\label{Lemma1.1}
		If $\M=g\mathcal{K}_u$ is a nearly invariant subspace with extremal function $g$ and associated
		inner function $u$ satisfying $u(0)=0$, then
		\begin{equation*}
			 P_{\M} f=g P_u(\overline{g} f)
		\end{equation*}
		for all $f$ in $\M$.  Consequently,
		the reproducing kernel for $\M$ is given by
		\begin{equation*}
			k_{\lambda}^{\M}(z) = \overline{g(\lambda)} g(z) \frac{1 - \overline{u(\lambda)} u(z)}{1 - \overline{\lambda} z}.
		\end{equation*}
	\end{lemma}

	Now armed with the preceding lemma, we are in a position to introduce truncated Toeplitz operators on nearly
	invariant subspaces. Certainly whenever $\phi$ is a bounded function we can use Lemma \ref{Lemma1.1}
	to see  that the operator $A_{\phi}^{\M}:\M\to\M$
	\begin{equation*}
		 A^{\M}_{\phi} f: =P_{\M}(\phi f)=gP_u(\overline{g}\phi f)
	\end{equation*}
	is well-defined and bounded.  More generally, we may consider symbols $\phi$ such that
	$|g|^2\phi$ belongs to $L^2$. In this case, for each $h$ in $\K_u^{\infty}:=\K_u\cap H^{\infty}$,
	the function $|g|^2\phi h$ is in $L^2$ whence $P_u(|g|^2\phi h)$ belongs to $\K_u$.
	By the isometric multiplier property of $g$ on $\mathcal{K}_u$, we see that
	\begin{equation*}
		P_{\M}(\phi h)  = g P_u(|g|^2\phi h)  \in g\K_u= \M.
	\end{equation*}
	Note that by the isometric property of $g$,
	the set $g\mathcal{K}_{u}^{\infty}$ is dense in $g\mathcal{K}_u$ by Theorem \ref{TheoremAleksandrovDensity}.
	Thus in this setting the operator  $A^{\M}_{\phi}$ is densely defined.  We refer to any such operator as a
	\emph{truncated Toeplitz operator} on $\M$.  We denote by
	$\mathscr{T}_{\M}$ the set of all such densely defined truncated Toeplitz operators
	which have bounded extensions to $\M$.
	The following theorem from \cite{RH}, which relies heavily upon the unitarity of the map \eqref{eq-IsometricMultiplier},
	furnishes the explicit link between $\mathscr{T}_{\M}$ and $\TTO_u$.

	\begin{theorem}[Hartmann-Ross]\label{prop1.2}
		If $\M=g\mathcal{K}_u$ is a nearly invariant subspace with extremal function $g$ and associated
		inner function $u$ satisfying $u(0)=0$, then for any Lebesgue measurable $\phi$ on $\dD$
		with $|g|^2 \phi \in L^2$ we have
		\begin{equation*}
			W_{g}^{*} A^{\M}_{\phi} W_{g} = A_{|g|^2 \phi}^{u}.
		\end{equation*}
	\end{theorem}

	In light of the preceding theorem, we see that the map
	\begin{equation*}
		A_{\phi}^{\M} \mapsto A_{|g|^2 \phi}^{u},
	\end{equation*}
	is a spatial isomorphism between $\mathscr{T}^{M}$ and $\TTO_u$.  In particular, we have
	\begin{equation*}
		\mathscr{T}^{M} = W_g \TTO_u W_g^*.
	\end{equation*}
	One can use the preceding results to prove the following facts about $\TTO^{M}$, all of which
	are direct analogues of the corresponding results on $\TTO_u$.
	\begin{enumerate}\addtolength{\itemsep}{0.5\baselineskip}
		\item $\TTO_{\M}$ is a weakly closed linear subspace of $\B(\M)$.
		\item $A^{\M}_{\phi} \equiv 0$ if and only if $|g|^2 \phi \in u H^2 + \overline{u H^2}$.
		\item  $C_{g} := W_g C W_{g}^{*}$ defines a conjugation on $\M$ and
			$A=C_{g} A^* C_{g}$ for every $A \in \TTO^{\M}$.
		\item If $S_{g} := W_g A_z W_{g}^{*}$, then a bounded operator $A$ on $\M$
			belongs to $\TTO^{\M}$ if and only if there are functions $\phi_1, \phi_2 \in \M$ so that
			\begin{equation*}
				A - S_{g} A S_{g}^{*} = (\phi_1 \otimes k^{\M}_{0}) + (k^{\M}_{0} \otimes \phi_2).
			\end{equation*}

		\item The rank-one operators in $\TTO^{\M}$ are constant multiples of
			\begin{equation*}
				 gk_{\lambda}\otimes gCk_{\lambda}, \qquad gCk_{\lambda}\otimes gk_{\lambda}, \qquad
				 gk_{\zeta}\otimes gk_{\zeta}.
			\end{equation*}

		\item $\TTO^{\M_1}$ is spatially isomorphic to $\TTO^{\M_2}$ if and only if either
			$u_1 = \psi \circ u_2 \circ \phi$ or $u_1 = \psi \circ \overline{u_2(\overline{z})} \circ \phi$
			for some disk automorphisms $\phi, \psi$.  In particular, this
			is completely independent of the corresponding extremal functions $g_1$ and $g_2$
			for $\M_1$ and $\M_2$.
	\end{enumerate}
\begin{acknowledgement}
First author partially supported by National Science Foundation Grant DMS-1001614.
\end{acknowledgement}

\def\cprime{$'$} \def\cprime{$'$} \def\cprime{$'$}
\providecommand{\bysame}{\leavevmode\hbox to3em{\hrulefill}\thinspace}
\providecommand{\MR}{\relax\ifhmode\unskip\space\fi MR }
\providecommand{\MRhref}[2]{%
  \href{http://www.ams.org/mathscinet-getitem?mr=#1}{#2}
}
\providecommand{\href}[2]{#2}

\end{document}